\documentclass{amsart}
\usepackage{amssymb}
\usepackage{amsfonts}
\usepackage{amsbsy}

\usepackage{latexsym}
\usepackage[mathscr]{eucal}
\usepackage{verbatim}
\usepackage{color}

\theoremstyle{plain}
\newtheorem{theorem}{Theorem}[section]

\newtheorem{corollary}[theorem]{Corollary}
\newtheorem{lemma}[theorem]{Lemma}
\newtheorem{proposition}[theorem]{Proposition}
\newtheorem{remark}[theorem]{Remark}
\newtheorem{definition}[theorem]{Definition}
\newtheorem{example}[theorem]{Example}

\newcommand{\A}{\mathcal{A}}
\newcommand{\B}{\mathcal{B}}
\newcommand{\C}{\mathcal{C}}

\newcommand{\K}{\mathcal{K}}
\newcommand{\h}{\mathcal{H}}
\newcommand{\J}{\mathcal{J}}
\newcommand{\I}{\mathcal{I}}
\newcommand{\eps}{\epsilon}

\DeclareMathOperator{\TA} {\mathcal{T}(\A)}
\DeclareMathOperator{\Aff} {Aff(\TA)}

\DeclareMathOperator{\LAff} {LAff(\TA)}

\DeclareMathOperator{\Ext} {\partial_{e}(\TA)}

\DeclareMathOperator{\her} {her}

\DeclareMathOperator{\tr} {Tr}

\DeclareMathOperator{\St} {\A\otimes \K}
\DeclareMathOperator{\M}{\mathcal M(\mathcal A \otimes \mathcal K)}
\newcommand{\Ma}{\mathcal M(\mathcal A)}
\newcommand{\Mb}{\mathcal M(\mathcal B)}
\DeclareMathOperator{\spn}{span}
\DeclareMathOperator{\rank}{rank}
\DeclareMathOperator{\Mul}{\mathcal M}
\DeclareMathOperator{\Ped}{Ped(\A)}

\newcommand{\Imin}{I_{\mathrm{min}}}
\newcommand{\Icon}{I_{\mathrm{cont}}}
\newcommand{\Io}{\mathcal{J}_o}
\DeclareMathOperator{\drr}{drr}
\DeclareMathOperator{\id}{id}

\def\sideremark#1{\ifvmode\leavevmode\fi\vadjust{\vbox to0pt{\vss
\hbox to 0pt{\hskip\hsize\hskip1em
\vbox{\hsize2cm\tiny\raggedright\pretolerance10000
\noindent#1\hfill}\hss}\vbox to8pt{\vfil}\vss}}}

\newcommand{\be}{\begin{equation}\label}
\newcommand{\ee}{\end{equation}}
\newcommand{\bq}{\begin{equation*}}
\newcommand{\eq}{\end{equation*}}
\newcommand{\ba}{\begin{align*}}
\newcommand{\ea}{\end{align*}}
\newcommand{\bp}{\begin{proof}}
\newcommand{\ep}{\end{proof}}
\newcommand{\bL}{\begin{lemma}\label}
\newcommand{\eL}{\end{lemma}}
\newcommand{\bP}{\begin{proposition}\label}
\newcommand{\eP}{\end{proposition}}
\newcommand{\bC}{\begin{corollary}\label}
\newcommand{\eC}{\end{corollary}}
\newcommand{\bT}{\begin{theorem}\label}
\newcommand{\eT}{\end{theorem}}
\newcommand{\bR}{\begin{remark}\label}
\newcommand{\eR}{\end{remark}}
\newcommand{\bD}{\begin{definition}\label}
\newcommand{\eD}{\end{definition}}
\newcommand{\bE}{\begin{example}\label}
\newcommand{\eE}{\end{example}}
\numberwithin{equation}{section}
\thanks{This work was partially supported by the Simons Foundation (grant No 245660 to Victor Kaftal and grant No  281966 to Shuang Zhang)}
\title[ ]{The minimal ideal in multiplier algebras}
\author{Victor Kaftal}

\address{Department of Mathematics\\
University of Cincinnati\\
P. O. Box 210025\\
Cincinnati, OH\\
45221-0025\\
USA}

\email{victor.kaftal@math.uc.edu}

\author{P.  W. Ng}

\address{Department of Mathematics\\
University of Louisiana\\
217 Maxim D. Doucet Hall\\
P.O. Box 41010\\
Lafayette, Louisiana\\
70504-1010\\
USA}

\email{png@louisiana.edu}

\author{Shuang Zhang}

\address{Department of Mathematics\\
University of Cincinnati\\
P.O. Box 210025\\
Cincinnati, OH\\
45221-0025\\
USA}

\email{zhangs@math.uc.edu}

\date{5/11/2017}


\begin{document}
\begin{abstract}
Let $\A$ be a simple,  $\sigma$-unital, non-unital, non-elementary C*-algebra and let  $\Imin$ be the intersection of all the ideals of $\Ma$ that properly contain $\A$. $\Imin$ coincides with the ideal  defined by Lin \cite{LinContScale} in terms of approximate units of $\A$ and $\Imin/\A$ is purely infinite and simple. If $\A$ is separable, or if $\A$ has the (SP) property and its dimension semigroup $D(\A)$ of Murray-von Neumann equivalence classes of projections of $\A$ is order separable, or if $\A$ has strict comparison of positive elements by traces, then $\A\ne \Imin$. 

If the tracial simplex $\TA$ is nonempty, let $\Icon$ be the closure of the linear span of the elements $A\in \Ma_+$ such that the evaluation map $\hat A(\tau)=\tau(A)$ is continuous. If $\A$ has strict comparison of positive element by traces then $\Imin=\Icon$. Furthermore, $\Imin$ too has strict comparison of positive elements in the sense that if $A, B\in (\Imin)_+$, $B\not \in \A$ and $d_\tau(A)< d_\tau(B)$ for all $\tau\in \TA$ for which $d_\tau(B)< \infty$, then $A\preceq B$.

However if $\A$ does not have strict comparison of positive elements by traces then $\Imin\ne \Icon$ can occur: a counterexample is provided by Villadsen's AH algebras without slow dimension growth. If the dimension growth is flat, $\Icon$ is the largest proper ideal of $\Ma$.

\end{abstract}

\maketitle  
\section{Introduction}
The ideal structure of the multiplier algebra of a simple, $\sigma$-unital non-unital non-elementary,  C*-algebra has received over the years a lot of attention. In this paper we will focus on the study of the smallest (closed) ideal properly containing $\A$.

 Lin \cite [Lemma 2]{LinIdealsAF} gave a constructive proof of the existence of such a smallest ideal for AF algebras in terms of the tracial simplex of the algebra (see \S 2.2). 
 
 Then  Lin and Zhang \cite{LinZhang},  proved that every simple, separable, non-unital, non-elementary C*-algebra with property (SP) and with an approximate identity of projections (such algebras do not need to have real rank zero) contains an $\ell^1$-sequence of projections (see Definition \ref {D:rapidly decreasing} for a generalization). Furthermore, all the principal ideals generated by projections associated to such sequences coincide with the minimal ideal properly containing $\A$.  
 
 In  \cite{LinContScale}  Lin defined for every simple $\sigma$-unital C*-algebra an ideal $\I$   in terms of an approximate identity of positive elements and proved that $\I$ is contained in any ideal properly containing $\A$.  If $\A$ is separable, then $A\ne \I$. This construction is instrumental in proving that $\Ma/\A$ is simple if and only if it is purely infinite and simple if and only if $\A$ has a continuous scale  (\cite [Theorems 2.4 and 3.2] {LinSimple}).

For simple C*-algebras with real rank zero, stable rank one, and weakly unperforated $K_0$, (equivalently, strictly unperforated monoid $V(\A)$ of Murray-von Neumann equivalence classes of projections in $\St$) Perera proved that there is a lattice  isomorphism between the ideals of $\Ma$ and the order ideals of $V(\A) \sqcup W_\sigma^d(S_u)$ (see \cite [Theorems 2.1 and 3.9] {PereraIdeals} and notations therein) and then proved  \cite[Proposition 4.1] {PereraIdeals}  that 
$V(\A) \sqcup  {\rm Aff}_{++}(S_u)$ is the smallest order ideal properly containing $V(\A)$, thus obtaining the smallest ideal properly containing $\A$. Here ${\rm Aff}_{++}(S_u)$ is the space of strictly positive  {\it continuous} affine functions on the state space $S_u$; see of also \S2.2 and \S 4.  This ideal, denoted by $L(\A)$, plays an important role in the study by Perera \cite{PereraIdeals} and  Kucerovsky and Perera  \cite{KuchPer} of the ideal structure of the multiplier algebra and the characterization of when the corona algebra $\Ma/\A$ is purely infinite. 

The goal of this paper is to clarify the relations between the various constructions of the minimal ideal and to further investigate its properties. Throughout the paper, $\A$ will denote a simple, $\sigma$-unital, non-unital and non-elementary C*-algebra. 

We revisit Lin's definition (\cite [Lemma 2.1]{LinContScale}) of a nonclosed ideal of $\Ma$ defined in terms of an approximate unit  $\{e_n\}$ of positive elements, which we denote by $L(K_o(\{e_n\}))$ and Lin did denote by $I_0$. It is easy to see that $L(K_o(\{e_n\}))$ is a left ideal.  Lemma 2.1 \cite {LinContScale} states that it is also *-invariant, and hence, a two-sided ideal, however Example \ref {E: Ko not invariant} shows a case when $L(K_o(\{e_n\}))$ is not two-sided. Nevertheless, its closure, that we  denote by $\Imin$ (and Lin did denote by $I$), is proved to be indeed a two-sided ideal (Corollary \ref {C: Imin}), and thus all the further results in \cite {LinContScale} and \cite {LinSimple} are correct. The proof that $\Imin$ is two-sided depends on a bidiagonal decomposition result (Theorem \ref {T:bidiag revisited}) which is in the line of the tri-diagonal decomposition of elements in $\Ma$, first introduced by Elliot in \cite [proof of Theorem 3.1]{ElliottAnn}.  More background on bi-diagonal and tri-diagonal decompositions is presented before Theorem \ref {T:bidiag revisited}. 
As a consequence of the proof, one also sees that $\Imin$ does not depend on the approximate identity chosen.  

In \cite [Remark 2.9]{LinContScale}) Lin did prove that $\Imin$ is contained in the intersection $\Io$ of all the ideals properly containing $\A$. In Theorem \ref {T: minimal} we prove that $\Imin= \Io$ and in  Theorem \ref {T:purely infinite} we show  that if $\Imin/\A$ is nonzero (and necessarily simple), then it is purely infinite.

Furthermore, $\A\ne \Imin$ if and only if there exists a  {\it thin} sequence of positive elements for $\A$ (Definition \ref {D:rapidly decreasing}, Theorem \ref {T:thin}). This notion can be seen as a generalization of the notion of $\ell^1$ sequence of projections introduced for the (SP) case in \cite{LinZhang}, thus providing a bridge between the approaches in \cite{LinZhang} and  \cite {LinContScale}.

If $\A$ is separable, or if $\A$ has the (SP) property and the dimension semigroup of Murray-von Neumann equivalence classes of projections is countable, or, more generally, is just order separable, then a thin sequence exists, and hence, $\A\ne \Imin$. This includes the case of type II$_1$ factors. 

We do not have examples when $\A=\Imin$. A natural test case we did consider is the nonseparable simple C*-algebra with both a nonzero finite and  an infinite projection studied by Rordam in \cite{RordamEx}, but it still yields $\A\ne \Imin$ (see last paragraph of Section \ref {S:strict comp}).

In the case when $\A$ has a nonempty tracial simplex $\TA$, another natural ideal inspired by the approaches in \cite {LinIdealsAF} and \cite {PereraIdeals} is $\Icon$, the ideal generated by positive elements with continuous evaluation function over $\TA$ (Definition \ref {D:Icon}). We show that  $\A\subsetneq \Icon$ (Proposition \ref {P:Icon not A}). If in addition, $\A$ has strict comparison of positive elements by traces, then $\Imin =\Icon$,  and hence, $\A\ne \Imin$ (Theorem \ref {T: Imin=Icont}). 
This result can be seen as a generalization of Perera's construction  \cite {PereraIdeals} of the minimal ideal in the case that  all quasitraces of $\A$ are traces (e.g., $\A$ is exact), while the weak unperforation of the $K_0$ group is equivalent to strict comparison by quasitraces, and hence, to strict comparison by traces.

What happens when there is no strict comparison by traces? In the case of the AH-algebras without slow dimension growth studied by Villadsen, which are known to have perforation, we prove that $\Imin \ne \Icon$ (Theorem \ref {T: Imin not Icont}). In addition, we show that if $\A$ has flat dimension growth, every positive element not in $\Icon$ must be full (Theorem \ref {T: Icon maximal}), and hence, $\Icon$ contains every other proper ideal of $\Ma$. If however the dimension growth is very fast, then this is no longer true (Proposition \ref {P: Rinfty}).

Finally, we prove that  if $\A$ has strict comparison of positive elements, then so does $\Imin$.  This result extends a previous result  obtained by us in the case when $\A$ is separable and has real rank zero (\cite [Proposition 3.1]{KNZCompJOT}). The methods used are inspired by the techniques used in \cite [Theorem 6.6]{KNZCompPos} to prove that $\Ma$ has strict comparison of positive elements if so does $\A$  and $\A$ has quasicontinuous scale in the sense of \cite{KuchPer}.

\section{Preliminaries}

\subsection{Cuntz subequivalence}
Cuntz subequivalence in a C*-algebra $\B$ is denoted by $\preceq$, that is, if $a, b\in \B_+$, then $a\preceq b$ if there is a sequence  $x_n\in \B$ such that $\|x_nbx_n^*-a\|\to 0$.  If $a\preceq b$ and $b\preceq a$, then $a$ is said to be equivalent to $b$ ($a\sim b$). It is well known that for projections subequivalence in this sense coincides with   Murray-von Neumann  subequivalence. 

We will use the following notation:
\be {e:feps}f_\epsilon(t):=\begin{cases}0&\text{for }t\in [0, \epsilon]\\
\frac{t-\eps}{\eps}&\text{for }t\in (\eps, 2\eps]\\
1&\text{for }t\in (2\epsilon, \infty).
\end{cases} \ee 

 For ease of reference we list here the following well known facts (see for instance \cite{CuntzDim},\cite {RordamStructure}).

\bL{L:Cuntz subeq}
Let $\B$ be a C*-algebra, $a,b\in \B_+$, $x\in \B$, $\delta>0$. Then
\item [(i)]  $xax^*\preceq a$;
\item [(ii)] $xx^*\sim x^*x$; 
\item [(iii)] If $a\le b$ then $a\preceq b$; 
\item[(iv)]   If  $\|a-b\|< \delta$ then $(a-\delta)_+\preceq b$; 
\item [(v)] If $a\preceq b$, then there is $r\in \B$ and $\delta'>0$ such that  $(a-\delta )_+= r(b-\delta')_+r^*$; there is also $r'\in \B$ such that $(a-\delta )_+= r'br'^*$
\item [(vi)]  If $a\preceq a'$ and $b\preceq b'$ then $a+b\preceq a'\oplus b'$;
\item [(vii)] If $ab=0$, then $(a+ b-\delta)_+=(a-\delta)_+ +(b-\delta)_+$;
\item [(viii)] \cite [Lemma 2.3] {KNZCompPos} If $a\le b$, then $(a-\delta)_+\preceq (b-\delta)_+$;
\item [(ix)] \cite [Lemma  2.4. (iii)]{KNZCompPos} $(a+b-\delta_1-\delta_2)_+\preceq(a-\delta_1)_++(b-\delta_2)_+$ for  $\delta_1, \delta_2\ge 0$.  
\eL

\bL{L: handy ineq}
Let $\B$ be a C*-algebra, $a, b\in \B_+$, and $\|a-b\|< \delta$. Then for all $\eps \ge 0$, $(a-\eps-\delta)_+\preceq (b-\eps)_+$.
\eL
\bp
Since $\|a-(b-\eps)_+\|\le \|a-b\|+\|b-(b-\eps)_+\| < \eps + \delta$, the conclusion follows from Lemma \ref {L:Cuntz subeq} (iv).
\ep
\bL{L: approx} Let $\B$ be a C*-algebra and $a\in  \B_+$. For every $\eps>0$ there is $y\in \B$ such that $\|a- a^{1/2}yay^* a^{1/2}\|< \eps$  and $\|yay^*\|=1$
\eL
\bp
Choose $g_\eps(t):=\sqrt{ \frac{ f_{\eps} (t) }{t} }$ and set  $y=g_\eps(a)$.  Then  $yay^* = f_{\eps}(a)$ and $$a- a^{1/2}yay^* a^{1/2}= a(1- f_{\eps}(a)),$$ hence both  conditions are satisfied.
\ep

We need the following results for which we do not have handy references. A related  result is \cite [Lemma 2.3] {LinContScale}.
\bL{L:smaller elements}
Let $\B$ be a simple C*-algebra and $0\ne a, b\in \B_+$. Then there is $0\ne c\in \B_+$ such that $c\preceq a$ and $c\le b$.
\eL

\bp
Since $\B$ is simple, there are elements $x_k, y_k\in \B$ such that
$$
\|\sum_{j=1}^nx_kay_k-b\|< \frac{\|b\|}{2}.
$$
Then $\sum_{j=1}^nx_kay_kb\ne 0$, and hence, there is some $k$ such that $x_kay_kb\ne 0$. Then also
$$c:= (x_kay_kb)^*(x_kay_kb)\ne 0, ~ d:= (x_kay_kb)(x_kay_kb)^*\ne 0.$$
First notice that $$d\le \|b\|^2\|y_k\|^2\|a\|x_kax_k^*\preceq  a, $$ whence $d\preceq a$. Since  $c\sim d$ by Lemma \ref {L:Cuntz subeq} (ii), it follows that  $c\preceq a$. On the other hand, 
$$c\le \|a\|^2\|x_k\|^2\|y_k\|^2\|b\|b,$$ hence $c\le b$, by scaling if  necessary $c$, which preserves the relation $c\preceq a$. 

\ep
For the convenience of the readers, we give the proof of  the following  well known results.
\bL{L:direct sums}
Let $\B$ be a simple, non-elementary C*-algebra. Then for every element $0\ne a\in \B_+$ there is an infinite sequence of mutually orthogonal elements $0\ne a_k\in \B_+$ such that $\sum_{j=1}^n a_k\le  a$ for all $n$.
\eL
\bp
Choose $\delta>0$ such that $(a-\delta)_+\ne 0$. Then $\her ((a-\delta)_+)$   contains a positive element $b$ with infinite spectrum (e.g., \cite [1.11.45]{LinAmenable}. In fact it contains an element with spectrum $[0,1]$ by \cite [pg 67] {AkemannSchultz}).  Since $b\le \frac{\|b\|}{\delta}a$, to simplify notations,  assume   that $b\le a$. Now choose by compactness a converging sequence of distinct elements in $t_k\in\sigma (b)$, and by passing to a subsequence assume that the $t_k$  are monotone and that $|t_k-t_{j+1}|$ is also monotone. Let $\epsilon_k:= \frac{1}{5}| t_k-t_{k+1}|.$ Then the intervals $[t_k-2\epsilon_k,t_k+2\epsilon_k]$ are disjoint.  Let $g_k$ be the continuous function with
$$g_k(t):=\begin{cases}0& t\in [0, t_k-2\epsilon_k]\cup[t_k+2\epsilon_k, \infty)\\
t_k-\epsilon_k& t\in [t_k-\epsilon_k,t_k+\epsilon_k]\\
\text{linear}&t\in [t_k-2\epsilon_k,t_k-\epsilon_k]\\
\text{linear}&t\in [t_k+\epsilon_k,t_k+2\epsilon_k]
\end{cases}
$$
 Let $a_k:=\frac{g_k(b)}{2^j}$. Then $0\ne a_k\le \frac{b}{2^j} \le \frac{a}{2^j}$ and $a_ia_k=0$ for $i\ne j$. Thus we conclude that $\sum_{j=1}^\infty a_k\le a$.
 \ep
 \bL{L:sqrt} Let $\B$ be a C*-algebra and let $a, b, c\in \B_+$ and $x\in \B$. Then
 \item [(i)] $xax^*\sim xa^2x^*$.
 \item [(ii)] $b^{1/2}ab^{1/2}\sim bab$.
 \item [(iii)] If $b\le c$, then $bab\preceq cac$. \eL
  \bp
 \item [(i)] First we see that $xa^2x^*\le \|a\|xax^*$ and hence $xa^2x^*\preceq xax^*$. For every $\delta>0$,  $0\le (a- \delta)_+\le \frac{1}{4\delta}a^2$ and hence $x(a- \delta)_+x^*\preceq xa^2x^*$. Thus $$xax^*= \lim_{\delta\to 0}x(a- \delta)_+x^*\preceq xa^2x^*,$$ which concludes the proof.
 \item [(ii)] \begin{alignat*}{2}b^{1/2}ab^{1/2}&\sim a^{1/2}ba^{1/2} &\text{(by Lemma \ref {L:Cuntz subeq} (ii))}\\
 &\sim a^{1/2}b^2a^{1/2} \text{(by  (i))}\\
 &\sim bab&\text{(by Lemma \ref {L:Cuntz subeq} (ii))}
 \end{alignat*}
\item [(iii)]
\begin{alignat*}{2}
bab&\sim b^{1/2}ab^{1/2}\qquad &(\text{by (ii)})\\
&\sim a^{1/2}ba^{1/2} \qquad &(\text{by Lemma \ref {L:Cuntz subeq} (ii)})\\
&\preceq a^{1/2}ca^{1/2} \qquad &(\text{by Lemma \ref {L:Cuntz subeq} (iii), since }a^{1/2}ba^{1/2}\le a^{1/2}ca^{1/2})\\
&\sim cac &(\text{by the same two equivalences above.})
\end{alignat*}
 \ep

\subsection{The tracial simplex and strict comparison}
Given a simple $\sigma$-unital (possibly unital) C*-algebra $\A$   and   a nonzero positive element $e$ in the Pedersen ideal $\Ped$ of $\A$, denote by $\TA$  the collection of the (norm) lower semicontinuous  densely defined tracial weights $\tau$ on $\A_+$, that are normalized on $e$. Explicitly, a trace $\tau$ is an additive and homogeneous map from $\A_+$   into $[0, \infty]$ (a weight); satisfies the trace condition $\tau(xx^*)=\tau(x^*x)$ for all $x\in \A$; the cone  $\{x\in \A_+\mid  \tau(x)< \infty\}$ is dense in $\A_+$  ($\tau$ is also called densely finite, or semifinite); satisfies the lower semicontinuity condition $\tau(x)\le \varliminf \tau(x_n)$ for $x, x_n\in \A_+$ and $\|x_n-x\|\to 0$, or equivalently,  $\tau(x)= \lim \tau(x_n)$ for $0\le  x_n\uparrow x$ in norm; and $\tau(e)=1$ ($\tau$ is normalized on $e$). 
We will assume henceforth that $\TA\ne \emptyset$, and hence, $\A$ is stably finite. 

When equipped with the topology of pointwise convergence on $\Ped$,  $\TA$ is a Choquet simplex (e.g., see \cite [Proposition 3.4]{TikuisisToms} and \cite {ElliottRobertSantiago}). The collection of the extreme points of $\TA$ is denoted by $\Ext$ and is called the {\it extremal boundary} of $\TA$. For simplicity's sake we call the elements of $\TA $ (resp., $\Ext$) traces (resp., extremal traces.) Tracial simplexes $\TA$ arising from different nonzero positive elements in $\Ped$ are homeomorphic; so we will not reference explicitly which element $e$ is used. A trace $\tau$ on $\A$ is naturally extended to the trace $\tau\otimes\tr $ on $\St$, and so we can identify  $\mathcal{T}(\St)$ with $\TA$. For more details, see  \cite {TikuisisToms}, \cite {ElliottRobertSantiago} and also  \cite {KNZCompPos} and  \cite {KNZComm}.

Recall also that as remarked in  \cite[5.3]{KNZComm}, by the work of F. Combes   \cite [Proposition 4.1, Proposition 4.4] {Combes} and Ortega, Rordam, and Thiel
\cite[Proposition 5.2]{OrtegaRordamThiel}
every $\tau \in \TA$ has a unique extension, (which we will still denote by  $ \tau$) to a lower semicontinuous (i.e., normal) tracial weight  (trace for short)  on the enveloping von Neumann algebra $\A^{**}$, and hence to a trace on the multiplier algebra $\Ma$.

\bD{D:Aff}
Given a convex compact space $K$,
\item [(i)]$Aff(K)$ denotes the Banach space of the continuous real-valued affine functions on $K$ with the uniform norm;
\item [(ii)]$LAff(K)$ denotes the collection of the lower semicontinuous affine functions on $K$ with values in $\mathbb R\cup\{+\infty\}$;
\item [(iii)] $Aff(K)_{++}$ (resp., $LAff(K)_{++}$) denotes the cone of the strictly positive functions (i.e., $f(x)>0$ for all $x \in K$) in $Aff(K)$ (resp., in $LAff(K)$.)
 \eD

For every $A\in \Ma_+$, denote by $\hat A$ the evaluation map
\be{e: eval map}\TA\ni\tau\to \hat A(\tau):=\tau(A)\in [0,\infty],\ee and denote by $\widehat{[A]}$ the dimension map
\be{e: dim map}\TA\ni\tau\to \widehat{[A]}(\tau):= d_\tau(A)\in [0,\infty]\ee
where $$ d_\tau(A):=\lim_n \tau(A^{1/n})$$ is the dimension function. 

Then it is well known that  $\hat A\in \LAff_{++}$ and $ \widehat{[A]}\in \LAff_{++}$ for every $A\ne 0$. By definition of the topology on $\TA$, if $a\in \Ped$, then $\hat a\in \Aff$.

As shown  in \cite[Remark 5.3]{OrtegaRordamThiel}, \be{e:range proj} d_\tau(A)=\tau(R_A)\quad\text{where }R_A\in \A^{**}\text{ is the range projection of }A.\ee 
We will also use frequently the following well known facts. If $A, B \in \Ma_+$, and $\tau\in \TA$ then
\begin{align}\label {e: ineq->ineq}  A\leq  B~ &\Rightarrow ~ \hat A(\tau) \le \hat B(\tau)\\
\label{e: suub->ineq}
A\preceq B~ &\Rightarrow ~ d_\tau(A)\le d_\tau(B) \\
AB=0~ &\Rightarrow ~ d_\tau(A+B)=d_\tau(A)+d_\tau(A)\label{e: dim sum}\\
\tau(A)&\le \|A\|d_\tau (A) \label{e:ineq 1}\\
 d_\tau((A-\delta)_+)&< \frac{1}{\delta}\tau(A)~\,~\forall  ~\delta >0. \label{e:ineq 2}
\end{align}

We will use the following notions of strict comparison.
\bD{D:str comp}
Let $\A$ be a simple C*-algebra with $\TA\ne \emptyset$. Then we say that 
\item [(i)] $\A$ has strict comparison of positive elements by traces if  
$a\preceq b$ for $a, b\in \A_+$ such that $d_\tau(a)< d_\tau(b)$ for all those $\tau\in \TA$ for which $d_\tau(b)< \infty$.
\item [(ii)] $\Ma$  has strict comparison of positive elements by traces if $A\preceq B$ whenever
$A, B\in \Ma_+$, $A$ belongs to the ideal $I(B)$ generated by $B$, and  $d_\tau(A)< d_\tau(B)$ for all those $\tau\in \TA$ for which $d_\tau(B)< \infty$. 
\eD
Notice that strict comparison is often defined in terms of 2-quasitraces. In \cite [Theorem 2.9] {KNZComm} we proved  that if a unital simple C*-algebra of real rank zero and stable rank one has strict comparison of positive elements by traces (equivalently, of projections, due to real rank zero)  then all 2-quasitraces are traces.  Recently this result was extended by showing that if a simple stable C*-algebra has strict comparison of positive elements by traces then all 2-quasitraces are traces \cite {NgRobert}.

Notice also that for $\Ma\ne \A$, which is not simple, we must add the obviously necessary hypothesis that $A\in I(B)$ as that condition does not follow in general from the comparison condition. Indeed if there is an element $B\in \A_+$ with  $d_\tau(B)= \infty$ for all $\tau\in \TA$ (and this  is certainly the case when  $\A$ is stable) then  the condition $d_\tau(A)< d_\tau(B)$ for all those $\tau\in \TA$ for which $d_\tau(B)< \infty$ is trivially satisfied for every $A\in \Ma_+$  and yet $A \not\preceq B$.

\subsection{Cones and ideals in C*-algebras}
Let $\B$ be a C*-algebra and  $K\subset \B_+$. Set
\begin{align}\label {e:L}L(K)&:=\{x\in \B\mid x^*x\in K\} \\
\label {e:L*L} L(K)^*L(K)&:=\big\{\sum_{j=1}^n x_j^*y_j\mid x_j,y_j\in L(K), ~n\in \mathbb N\big\}.
\end{align}
\bD{D: cones} Let $\B$ be a C*-algebra and  $K\subset \B_+$.
\item [(i)] $K$ is a {\it  cone} if   $x+y\in K$ and $tx\in K$ whenever $x, y\in K$ and $0\le t\in \mathbb R$; $K$ is {\it hereditary}  if $x\in K$ whenever $0\le x\le y\in K$. 
\item [(ii)] A  subalgebra $\C\subset \B$ is hereditary if the cone $\C_+$ is hereditary.
\item [(iii)] A cone  $K$ is 
\begin {enumerate} \item [(a)] {\it invariant}  if $axa^*\in K$ whenever $x\in K$ and $a\in \B$;
\item [(b)] {\it strongly invariant} if  $x^*x\in K$ whenever $x\in \B $ and $xx^*\in K$;
\item [(c)] {\it weakly invariant} if $axa^*\in \bar K$ whenever $x\in K$ and $a\in \B$.
\end{enumerate}
\eD
Hereditary cones are also called {\it order ideals}. It is well known and immediate to see that if  $K$ is a hereditary cone, then $L(K)$ is a left ideal of $\B$, $L(K)^*L(K)$  and $L(K)^*\cap L(K)$ are *-subalgebras of $\B$, and $L(K)^*L(K)\subset L(K)^*\cap L(K)$.   Furthermore, if $K$ is a hereditary cone, then
\be{e:invar}  \text{$L(K)$ is two-sided if and only if $K$ is invariant.}\ee
\be{e:str invar}\text{ $L(K)= L(K)^*$  if and only if $K$ is strongly invariant.}\ee

\bT{T cones} Let $\B$ be a C*-algebra and  $K\subset \B_+$ be a hereditary cone. Then
\item [(i)] $\bar K$ is a hereditary cone (\cite [Theorem 2.5]{Effros63}) 
\item [(ii)] $L(K)^*L(K)=\text{span}\, K$ (the {\it the collection of complex  linear combinations of $K$}) and $(L(K)^*L(K))_+=K$  (\cite [Proposition 3.21] {StratilaZsido}.
  \item [(iii)] If $K$ is closed, then $L(K)^*L(K)= L(K)^*\cap L(K)$ and the mappings  $\B\to\B_+$, $K\to L(K)$, and $L\to L^*\cap L$ define bijective, order preserving correspondences between the sets of hereditary C*-subalgebras of $\B$, closed hereditary cones of $\B_+$, and closed left ideals of $\B$  (\cite [Theorem 2.4] {Effros63}, \cite [Theorem 1.5.2]{Pedersen}).
\eT

We collect here some simple properties of hereditary cones in C*-algebras that we will use in this paper.
\bL{L:new stuff cone} Let $\B$ be a C*-algebra and $K\subset \B_+$ be a cone.
\item [(i)] The (norm) closure $\bar K$ of $K$ is a cone.
\item [(ii)] If $K$ is weakly invariant, then $\bar K$ is  invariant. 
\item [(iii)] If $K$ is  invariant, then  $\bar K=\{x\in \B_+\mid (x-\delta)_+\in K ~\forall ~\delta>0\}.$
\item [(iv)] If $K$ is closed and invariant, then it is hereditary and strongly invariant.

\eL
\bp
\item [(i)] Obvious.
\item [(ii)]  Let $x\in \bar K$, $ a\in \B$ and let $x_n\in K$ be a sequence converging (in norm) to $x$. Since $ax_na^*\in \bar K$ for every $n$, it follows that $axa^*=\lim_n ax_na^*\in \bar K$, that is, $\bar K$ is invariant.
\item [(iii)] Let $K':=\{x\in \B_+\mid (x-\delta)_+\in K ~\forall ~\delta>0\}$. Since $\underset{\delta\to 0}{\lim}(x-\delta)_+= x$ for all $x\in \B_+$, it follows that $K'\subset \bar K$.
Conversely, let $x\in \bar K$, $\delta >0$, and choose $y\in K$ such that $\|x-y\|< \frac{\delta}{2}$. Then  $(x-\delta)_+=ryr^*\in K$ for some $r\in \B$ by Lemma \ref {L:Cuntz subeq} (iv) and (v). Thus  $x\in K'$, which proves that $K'=\bar K$.
\item [(iv)]  Let  $x\le y$, with  $x\in \B_+$ and $y\in K$.  By Lemma \ref {L:Cuntz subeq}(iii) and (v)  for every $\delta >0$ there is an $r\in \B$ such that $(x-\delta)_+=ryr^*\in K$ (because $K$ is invariant). Thus $x= \underset{\delta\to 0}{\lim}(x-\delta)_+ \in K$ (because $K$ is closed),  which proves that  $K$ is hereditary.  

Now let $x^*x\in K$ and $x= v|x|$ be the polar decomposition of $x$. By \cite [Lemma  2.1] {AkemannPedersen77}, $v|x|^{1/n}\in \B$ for every $n\in \mathbb N$, hence $(v|x|^{1/n})x^*x (v|x|^{1/n})^*\in K$ . Since $|x|^{1/n}|x|\to |x|$ in norm, it follows  that also  $$xx^*= vx^*xv^*=\lim_n(v|x|^{1/n})x^*x (v|x|^{1/n})^*\in K,$$ which proves that $K$ is strongly invariant.
\ep 

In the course of the proof of (iv) we have shown that
\be{e: Kinvar not her} \text{If $K$ is  invariant and $0\le x\preceq y\in K$ then $(x-\delta)_+\in K$ for all $\delta >0$.}\ee

From Example \ref {E: Ko not invariant} and Corollary \ref {C: Imin},   we will see that the condition in (iii) that $K$ is invariant cannot be replaced by condition that $K$ is weakly invariant. 

\bC{C:closure ideal}
Let $\B$ be a C*-algebra and   $K\subset \B_+$ is a weakly invariant hereditary cone in a C*-algebra $\B$,  then  $\overline {L(K)}=L(\bar K)$, $\overline {L(K)}$ is a two-sided ideal, and $\overline {L(K)}_+=\bar K$.
\eC

\bp

By Lemma \ref {L:new stuff cone} (i), (ii),  and (iv), $\bar K$ is a strongly invariant hereditary cone. By (\ref {e:str invar}), $L(\bar K)=L(\bar K)^*$ and by Theorem \ref {T cones} (i) and (ii), $\spn \bar K =L(\bar K)$.
Since $K$ is hereditary, $L(K)$ is a left ideal, and hence, so is $ \overline{L(K)}$. Moreover, $K\subset L(K)$, hence $\spn \bar K\subset \overline{L(K)}$, and hence, $L(\bar K)\subset \overline{L(K)}$. 
On the other hand, $L(K)\subset L(\bar K)$, and hence $ \overline{L(K)}\subset {L(\bar K)}$, and thus $ \overline{L(K)}= {L(\bar K)}$. 
 \ep

\subsection{Approximate identities}\label{S:approx id}
When $\B$ is a $\sigma$-unital C*-algebra, and $\{e_n\}$ is an approximate identity, we will always assume  that 
\be{e: appr id} \text{$\{e_n\}$ is strictly increasing ($0\le e_{n}\lneqq e_{n+1}$) and that } ~ e_{n+1}e_n=e_n\quad\forall\, n.\ee
  It is also convenient to define $e_0=0$. Notice that $e_n\in \Ped$ and $\|e_n\|=1$ for all $n\ge 1$.

Notice that 
\be{e:appr id 1}(e_{m+1}-e_{n-1})(e_m-e_n)= e_m-e_n\quad \forall ~m>n,\ee
and hence,
\be{e:appr id 2} e_m-e_n\le R_{e_m-e_n}\le e_{m+1}-e_{n-1} \quad \forall ~m>n. \ee
\bR{R: special approx id} We can always pass from  an approximate identity satisfying the above conditions  to a subsequence $f_n$ satisfying the following two stronger conditions assumed in \cite {LinContScale}:
\begin{itemize} \item [(a)] 
Let $g_n:=f_n-f_{n-1} $ (with $f_0:=0$), then $\|g_n\|=1$ for all $n$ and $g_ng_m=0$ for $|m-n|\ge 2$.
\item [(b)] There are $a_n\in \B_+$ with $\|a_n\|=1$ such that $a_n \le g_n,\, a_ng_n=g_na_n=a_n$ and $a_ng_m =0$ for $n\ne m.$
\end{itemize}
\eR
\bp
Let $f_n:= e_{5n}$. Clearly, $(f_n-f_{n-1})(f_m-f_{m-1})=0$ for $|m-n|\ge 2$. Set $a_n:= e_{5n-1}-e_{5n-4}$. Then $f_n-f_{n-1}\ge a_n$ by the monotonicity of $e_n$ and 
$$(f_n-f_{n-1})a_n=a_n(f_n-f_{n-1})=a_n$$
by (\ref{e:appr id 1}). Furthermore,  $\|a_n\|=1$ since by (\ref{e:appr id 2}), $a_n\ge R_{ e_{5n-2}-e_{5n-3}}\ne 0$; in particular, $\|f_n-f_{n-1}\|=1$. \ep

\section{The minimal ideal and its hereditary cone} 
\bD{D:cone Ko} \cite [Lemma 2.1] {LinIdealsAF} Let $\A$ be a simple, $\sigma$-unital, non-unital C*-algebra with an approximate identity $\{e_n\}$. Then we define the following set of positive elements in $ \Ma$:
 \begin{multline}\notag K_o(\{e_n\}):=\{X\in \Ma_+\mid \forall~0\ne a\in \A_+ ~\exists~N\in \mathbb N \\ \ni ~m>n \ge N\Rightarrow (e_m-e_n)X(e_m-e_n)\preceq a\}.\end{multline}  
\eD
\bR{R:equiv K_o} \item [(i)] By Lemma \ref {L:sqrt} (iii), 
\begin{multline}\notag K_o(\{e_n\}):=\{X\in \Ma_+\mid \forall~0\ne a\in \A_+ ~\exists~N\in \mathbb N \\ \ni ~m> N\Rightarrow (e_m-e_N)X(e_m-e_N)\preceq a\}.\end{multline}  
This equivalent  formulation  will also be  used in the paper. 
\item[(ii)] If $\A$ has the (SP) property, (i.e., every nonzero hereditary subalgebra of $\A$ contains a nonzero projection), then for every $0\ne a\in \A_+$ there is a projection $0\ne p\preceq a$. Thus in the defining property of $K_o(\{e_n\})$ we can replace ``for all  nonzero elements $a\in \A_+$" with ``for all  nonzero projections $p\in \A.$"
\eR

\bL{L:Ko hered cone}
\item [(i)] $X\in K_o(\{e_n\})$ if and only if $X^{1/2} \in K_o(\{e_n\})$.
 \item [(ii)]  $K_o(\{e_n\})$ is a hereditary cone of $\Ma$ if and only if $\A$ is non-elementary.

\eL
\bp
\item  [(i)] Immediate from the definition and Lemma \ref {L:sqrt} (i).
\item [(ii)] 
It is also immediate to verify that  $K_o(\{e_n\})$ is always hereditary and that if $X\in K_o(\{e_n\})$ then $tX\in K_o(\{e_n\})$ for every $t\ge 0$. Assume first that $\A$ is non-elementary and that  $X, Y \in  K_o(\{e_n\})$.  Let $0\ne a\in \A_+$, then by Lemma \ref {L:direct sums} we can find two elements $0\ne a', a''\in \A_+$  with $a'a''=0$ and $a'+a''\le a$. Let $N'$ (resp., $N''$) be such that for all $m>n\ge N'$ (resp.,  $m>n\ge N''$), we have $(e_m-e_n)X(e_m-e_n)\preceq a'$ (resp.,   $(e_m-e_n)Y(e_m-e_n)\preceq a''$).  Hence, for all $m>n\ge N:= \max (N', N'')$ we have by Lemma \ref {L:Cuntz subeq}(vi)
$$
(e_m-e_n)(X+Y)(e_m-e_n)= (e_m-e_n)X(e_m-e_n)+(e_m-e_n)Y(e_m-e_n)\preceq a'+a''\le a.
$$
Thus $K_o(\{e_n\})$ is a cone.

Assume now that $\A= \K$, and hence, $\Ma=B(\h)$, and let $\{e_n\}$ be an increasing sequence of rank $n$ projections. Then  it is easy to verify that
$$K_o(\{e_n\})= \{x\in B(\h)_+\mid \exists n\,  \ni\, \text{rank}(1-e_n)x(1-e_n)\le 1\}.$$ Let $\{\eta_n\}$ be an  orthonormal basis of $\h$ such that  span$ \{\eta_1, \cdots, \eta_n\}=R_{e_n}$, and let $\xi:= \sum_{j=1}^\infty \frac{1}{2^j}\eta_{2j}$ and $\xi':= \sum_{j=1}^\infty \frac{1}{2^j}\eta_{2j+1}$. Then both $\xi\otimes \xi$ and $\xi'\otimes \xi'$ belong to $K_o(\{e_n\})$ since they have rank one, but $(1-e_n)\big(\xi\otimes \xi+\xi'\otimes \xi'\big) (1-e_n)$ has rank two for every $n$, and hence, $\xi\otimes \xi+\xi'\otimes \xi'\not \in K_o(\{e_n\})$.
\ep

\bC{L:left ideal}
Let $\A$ be a simple, $\sigma$-unital, non-unital, non-elementary C*-algebra with an approximate identity $\{e_n\}$. Then $L(K_o(\{e_n\})$ is a left ideal and $$L(K_o(\{e_n\}))_+=K_o(\{e_n\}).$$
\eC
That $L(K_o(\{e_n\})$ is a left ideal follows immediately from the fact that $K_o(\{e_n\})$ is a  hereditary cone. The equality $K_o(\{e_n\}) = L(K_o(\{e_n\}))_+$ was suggested by H. Lin (private communications).
$L(K_o(\{e_n\})$ was denoted by $\mathcal I_o $ in \cite {LinContScale}. Contrary to what was stated in \cite [Lemma 2.1] {LinContScale}), the
 following example shows that $K_o(\{e_n\})$  is in general not invariant, i.e., the ideal $L(K_o(\{e_n\}))$ is not two-sided.

\bE{E: Ko not invariant} Let $\A_o$  be a simple, unital, finite non-elementary C*-algebra and let $\A:= \A_o\otimes \K$.  Let $\{e_{ij}\}$ be the standard matrix units in $\K$, then  $e_n:=1\otimes \sum_{k=1}^ne_{kk}$ is an  increasing approximate identity of projections of $\A$. Let $$V:= 1\otimes\sum_{k=1}^\infty 2^{-k/2} e_{1,k}.$$ Then $VV^*= e_1=1\otimes e_{11}\in  K_o(\{e_n\})$, i.e., $V^*\in L( K_o(\{e_n\}))$.
Let $$P:=V^*V= 1\otimes \sum_{h,k=1}^\infty 2^{-(h+k)/2} e_{h,k}.$$
For every $n>1$ and $0\ne a\in (\A_o)_+$ with $a\not\sim 1$ we have
$$(e_n-e_{n-1})P(e_n-e_{n-1})= 1\otimes 2^{-n} e_{n,n} \sim  1\otimes e_{11}\not\preceq a\otimes e_{11}.$$
Thus $P\not \in K_o(\{e_n\})$, i.e., $V\not \in L(K_o(\{e_n\}))$. This example shows that the cone $K_o(\{e_n\})$ is not invariant, and,  equivalently, that $L(K_o(\{e_n\}))$ is not a two-sided ideal. 
It also shows that $K_o(\{e_n\})$ does not satisfy the conclusion of Lemma \ref {L:new stuff cone} (iii) since $P \in \overline{K_o(\{e_n\})}$ and yet $\frac{1}{2}P=(P-\frac{1}{2})_+\not \in K_o(\{e_n\})$. Furthermore, if we choose an approximate identity $f_n= 1\otimes \sum_{k=1}^nf_{kk}$ with $f_{1,1}= \sum_{h,k=1}^\infty 2^{-(h+k)/2} e_{h,k}$, we see that $P\in K_o(\{f_n\})$, which shows that $K_o(\{e_n\})\ne K_o(\{f_n\})$.
\eE

In Corollary \ref {C: Imin}  we will see that $K_o(\{e_n\})$ is always weakly invariant, and hence, $\overline{K_o(\{e_n\})}$ is strongly invariant and that $\overline{K_o(\{e_n\})}$ does not depend on the approximate identity $\{e_n\}$.   Meanwhile, the next lemma shows that refinements of an approximate identity do not change the cone $K_o$.
 \bL{L:Ko for subseq}
Let $\A$ be a simple $\sigma$-unital non-unital non-elementary C*-algebra with  an approximate identity $\{e_n\}$. Then 
$K_o(\{e_n\})=K_o(\{e_{n_k}\})$  for any strictly increasing sequence  $n_k$  of integers.
 \eL
\bp
Let $X\in K_o(\{e_{n_k}\})$ and  $0\ne a\in \A_+$. Then there is an  $L\in \mathbb N$ such that if $k> L$ then
$(e_{n_k}-e_{n_L})X(e_{n_k}-e_{n_L})\preceq a$.  Let $m> n_L$ and choose $k$ such that $n_k\ge m$. Then $e_{m}-e_{n_L}\le e_{n_k}-e_{n_L}$, and hence, by Lemma \ref {L:sqrt} (iii)
$$(e_{m}-e_{n_L})X(e_{m}-e_{n_L})\preceq (e_{n_k}-e_{n_L})X(e_{n_k}-e_{n_L})\preceq a.$$
Thus $X\in K_o(\{e_n\})$.  The opposite inclusion is obvious.

\ep  

Given any approximate identity $\{e_n\}$ of $\A$, it is clear that $e_nae_n\in K_o(\{e_n\})$ for every  $a\in \A_+$ and $n\in \mathbb N$.  Since $e_nae_n\to a$, it follows that 
\be{e: Asubclos K0}\A_+\subset \overline{K_o(\{e_n\})}.\ee
The inclusion $\A_+\subset K_o(\{e_n\})$ is however equivalent to the condition that $\A$ has continuous scale.
Recall that $\A$ is said to have  continuous scale if for some (and hence,  for every) approximate identity $\{e_n\}$ and for every $0\ne a \in \A_+$ there is an $N\in \mathbb N$ such that $e_m-e_n\preceq a$ for all $m>n\ge N$. 

\bL{L:cont scale} Let $\A$ be a simple, $\sigma$-unital, non-unital, and non-elementary C*-algebra with  an approximate identity $\{e_n\}$. The following are equivalent.
\item[(i)] $\A$ has continuous scale;
\item [(ii)]  $K_o(\{e_n\})=\Ma_+$;
\item [(iii)] 	$\overline{K_o(\{e_n\})}=\Ma_+$;
\item [(iv)] $\A_+\subset K_o(\{e_n\}).$
\eL
\bp
(i) $\Rightarrow$ (ii). For every $x\in  \Ma_+$ and every $m>n$ we have
$$(e_m-e_n)x(e_m-e_n)\le \|x\| (e_m-e_n)\preceq e_m-e_n.$$
(ii) $\Rightarrow$ (iii) and (ii) $\Rightarrow$ (iv). Obvious\\
(iii) $\Rightarrow$ (ii). Since $1\in  \overline{K_o(\{e_n\})}$, there is an $x\in K_o(\{e_n\})$ such that $\|x-1\|< 1$. Thus $x$ is invertible, and hence, $\rho 1\le x$ for some scalar $\rho>0$. Since $ K_o(\{e_n\})$ is a hereditary cone, it follows that $1\in K_o(\{e_n\})$, hence $\Ma_+\subset K_o(\{e_n\})$ and thus  (ii) holds.\\
(iv) $\Rightarrow$ (i).
Let $b:= \sum _{k=1}^\infty \frac{1}{k} (e_{k+1}-e_k)$ where the convergence is in norm, and hence, $b\in \A_+\subset K_o$. Then for every $0\ne a \in \A_+$ there is an $N\in \mathbb N$ such that if $m\ge n\ge N+1$ then $(e_{m+1}-e_{n-1})b(e_{m+1}-e_{n-1})\preceq a$. But then by (\ref {e:appr id 1}) we also have
\ba
e_m-e_n&\sim \frac{1}{m-1}(e_m-e_n)\\
&= \frac{1}{m-1}(e_{m+1}-e_{n-1})(e_m-e_n)(e_{m+1}-e_{n-1})\\
&=(e_{m+1}-e_{n-1})\sum_{k=n}^{m-1}\frac{1}{m-1}(e_{k+1}-e_k)(e_{m+1}-e_{n-1})\\
&\le (e_{m+1}-e_{n-1})\sum_{k=n}^{m-1}\frac{1}{k}(e_{k+1}-e_k)(e_{m+1}-e_{n-1})\\
&\le (e_{m+1}-e_{n-1})b(e_{m+1}-e_{n-1})\\
&\preceq a.
\end{align*}
Thus the scale is continuous.
\ep
The implication (iii) $\Rightarrow$ (ii) is essentially the ``only if" part of \cite [Theorem 2.10] {LinContScale}. 
The following notions have appeared in various forms and various names in the literature (e.g., \cite{LinZhang}, \cite [4.3.11]{BlackB88}) and for ease of reference we present them by the following formal definition. 
\bD{D:rapidly decreasing} Let $\B$ be a C*-algebra. 
\item [(i)]  A sequence $0\ne s_i\in \B_+$ is called {\it order dense}  for $\B$ if for every $0\ne a\in \B_+$ there is an integer $n$ for which $s_n\preceq a$.
\item [(ii)] A sequence of mutually orthogonal elements $0\ne t_i\in \B_+$ is called thin for $\B$ if for every $0\ne a\in \B_+$ there is an integer $N$ such that $\sum_{i=n}^m t_i\preceq a$ for all $m\ge n\ge N$.
\eD
Recall that a thin sequence of projections is called an $\ell^1$ sequence in \cite {LinZhang}. Clearly, thin sequences are  order dense; also if $\{s'_i\}$ is an order dense sequence for $\B$ and $0\ne s_i\in \B_+$ with  $s_i\preceq s'_i$ for every $i$,  then $\{s_i\}$  is also order dense for $\B$. Similarly, 
let $0\ne s_i,  s_i'\in \B_+$.
\be{e: subthin}
\text {If $\{s'_i\}$ is thin, $ s_is_j=0$ for $i\ne j$ and $s_i\preceq s'_i~\, \forall \, i$, then $\{s_i\}$  is  thin.}
\ee
This follows from Lemma \ref {L:Cuntz subeq} (vi) since  $\sum_{i=n}^ms_i\preceq \sum_{i=n}^ms'_i$ for every $m\ge n$. It is also immediate to see that
\be{e: subthin2}
\text {If $\{s'_i\}$ is thin, $s_i= \alpha_i  s'_i$  for some $\alpha_i>0$, then  $\{s_i\}$  is  thin.}
\ee

 In separable C*-algebras, it easy to construct order dense sequences (see also the construction in \cite [Lemma 2.4]{LinContScale} and \cite {SZ90} for projections).
\bP{P:separable} Every separable C*-algebra has an order dense sequence.
\eP
\bp
Let $\B $ be a separable C*-algebra and let  $\{b_m\}$ be a sequence of  positive elements  dense in the unit ball of $\B_+$.   Let $\{s_n\}$ be an enumeration of the nonzero elements in the collection $\{(b_m-\frac{1}{2})_+\mid m\in \mathbb N\}$. For every $0\ne a\in \B_+$ there is an $m\in \mathbb N$  such that $\|\frac{a}{\|a\|}-b_m\|< \frac{1}{2}$. Then $\|b_m\|> \frac{1}{2}$, hence $(b_m-\frac{1}{2})_+\ne 0$, and thus $(b_m-\frac{1}{2})_+=s_n$ for some $n$. Then $s_n\preceq \frac{a}{\|a\|}\sim a$ by Lemma \ref {L:Cuntz subeq} (iv). Thus $\{s_n\}$ is an order dense sequence.

\ep
Another case when order dense sequences are immediate to obtain is the following. For every C*-algebra  $\A$, denote by $D(\A)$ the (possibly empty) dimension semigroup of Murray-von Neumann equivalence classes of projections. We say that $D(\A)$ is order separable if there is a sequence $p_n $ of nonzero projections of $\A$ such that for every projection $0\ne p\in \A$ there is a $p_n\preceq p$. Of course, if $D(\A)$ is countable, it is also order separable, but type II$_1$ von Neumann factors are examples of (non-separable) C*-algebras with a dimension semigroup $D(\A)$ that is order separable but not countable. 
\bP{P: countable dim group}
Every C*-algebra  $\B$  with (SP) property and with order separable dimension semigroup $D(\B)$  has an order dense sequence of projections. 
\eP
\bp
By the (SP) property, for every $0\ne a\in \B_+$ there is a nonzero projection $q\in \her(a)$, and hence, $q\preceq a$.  Since $p_n \preceq q$ for some $n$, we have $p_n \preceq a$. Thus $\{p_n\}$ is order dense for $\B$.  
\ep

The following construction permits to construct thin sequences starting with order dense sequences. For future use in this paper, we will prove a slightly stronger version than needed in this section.
When $s, t\in \B_+$ and $n\in \mathbb N$, we will denote by $ns$ an n-fold direct sum of $s$ with itself. 
Then $ns\in M_n(\B_+)$ and the subequivalence relation $ns\preceq t$ is understood to hold in $M_n(\B_+)$.
In particular, if $s\preceq t_i$ for $1\le i\le n$ and  $t_i$ are mutually orthogonal, then by Lemma \ref {L:Cuntz subeq} (vi), 
$ns\preceq \sum_{i=1}^n t_i$.
 \bL{L:ncopies}
 Let $\B$ be a simple non-elementary C*-algebra. Then for every sequence $s_i$ of elements $0\ne s_i\in \B_+$, there is a sequence of mutually orthogonal elements $0\ne t_i\in \B_+$ such that $n\sum_{i=n}^mt_i \preceq s_n$ for every pair of integers $m\ge n$.  \eL
\bp
Let $0\ne a_i\in \B_+$ be a sequence of mutually orthogonal elements (e.g., see Lemma \ref {L:direct sums}). By Lemma \ref {L:smaller elements}, there are elements $0\ne s'_i\in \B_+$ with $s'_i\le a_i$ and $s'_i\preceq s_i$. 
For every $i$, use Lemma \ref {L:direct sums} to find an infinite sequence of mutually orthogonal nonzero elements $0\ne s'_{i,j}\in \B$ such that $\sum_{j=1}^ns'_{i,j} \le s'_i$ for every $n$.  For every $j$, set $s_{1,j}= s'_{1,j}$.   Applying  Lemma \ref {L:smaller elements},   find an element   $0\ne s_{2,j}\le s'_{2,j}$, such that $s_{2,j}\preceq s_{1,j}$. By iterating the construction, find  sequences $0\ne s_{i,j}\le s'_{i,j}$ such that 
$$s_{i,j}\preceq  s_{i-1,j}\preceq \cdots\preceq s_{1,j}\quad \forall \, i,\, j.$$
Now apply again Lemma \ref {L:direct sums} to find mutually orthogonal elements $0\ne t_{i,j}\in \A_+$ such that $\sum_{j=1}^nt_{i,j}\le s_{i,i}$. By Lemma \ref {L:smaller elements} we can assume again that for every $i$.
$$t_{i,i}\preceq t_{i,i-1}\preceq \cdots \preceq t_{i,1}.$$
Let $t_i:=t_{i,i}$. Notice that the sequences $t_i\le s_{i,i}\le s'_i\le a_i$ are mutually orthogonal. Thus for every $n\le i\in \mathbb N$
$$nt_i\preceq \sum_{j=1}^{i} t_{i,j} \le s_{i,i},$$
and hence,
$$n\sum_{i=n}^mt_i \preceq \sum_{i=n}^m s_{i,i} \preceq \sum_{i=n}^m s_{n,i}\preceq \sum_{i=n}^m s'_{n,i}\le s'_n\preceq s_n.$$
\ep 
 
The following consequence is now immediate. 
 \bC{C: infinit vs thin} Let $\B$ be a simple non-elementary C*-algebra. 
If $\B$ has an order dense sequence $s_i$, then it has a thin sequence $t_i$ with $t_i\preceq s_i$ for every $i$. \eC

If $\A$ is a simple, $\sigma$-unital, non-unital C*-algebra with an approximate identity $\{e_n\}$, and two sequences of positive integers $m_j$ and $n_j$ such that $n_j< m_j< n_{j+1}$ for every $j$, $d_j\in \A_+$  is a bounded sequence for which $d_j\le M_j (e_{m_j}-e_{n_j})$ for some $M_j>0$, then the elements $d_j$ are mutually orthogonal and  the series $\sum_{j=1}^\infty d_j$ converges strictly. We will call the sum $D$ of such a sequence {\it diagonal with respect to $\{e_n\}$}. Furthermore, $D\in \A$ if and only if $\lim \|d_j\|=0$.

\bL{L: diag}
Let $\A$ be a simple, $\sigma$-unital, non-unital C*-algebra and assume that $D:=\sum_{j=1}^\infty d_j$ is diagonal with respect to an approximate identity $\{e_n\}$.  Then $D\in K_o(\{e_n\})$ if and only if the sequence $\{d_j\}$ is thin.
\eL
\bp
Let $n_j$ and  $m_j$ be sequences of positive integers such that $n_j< m_j< n_{j+1}$ for every $j$ and $d_j\le M_j (e_{m_j}-e_{n_j})$ for some $M_j>0$ and all $j$. 
Assume first that the sequence $\{d_j\}$ is thin.  Since for every $p\ge L\in \mathbb N$ we have 
$$(1-e_{n_L})\sum_{j=1}^{L-1}d_j=0\quad\text{and}\quad e_{m_p} \sum_{j=p+1}^{\infty}d_j=0$$
it follows that
$$(e_{m_p}-e_{n_L}) D= (e_{m_p}-e_{n_L}) \sum_{j=L}^{p}d_j.$$
Now let $0\ne a \in \A_+$ and $L\in \mathbb N$ be such that if $p\ge L$, then $ \sum_{j=L}^{p}d_j\preceq a$. For every $m>N:=n_L$ choose $p$ such that $m_p\ge m$. Then by Lemma \ref {L:sqrt} (iii) and Lemma \ref {L:Cuntz subeq} (i), we have
\ba (e_m-e_{N})D(e_m-e_{N})&\preceq (e_{m_p}-e_{N})D(e_{m_p}-e_{N})\\
&=  (e_{m_p}-e_{N})\sum_{j=L}^{p}d_j(e_{m_p}-e_{N})\\
&\preceq \sum_{j=L}^{p}d_j\\
&\preceq	 a.\end{align*}
This proves that $D\in K_o(\{e_n\})$.

Assume now that $D\in K_o(\{e_n\})$ and let $0\ne a\in \A_+$. Then there is an integer $N$ such that $(e_m-e_{N})D (e_m-e_{N})\preceq a$ for every $m\ge N$.
Let $L$ be such that $n_L\ge N+1$ and $p\ge L$.  Since $\sum_{j=L}^{p}d_j\le M(e_{m_p}-e_{n_L})$ and  $$(e_{m_p+1}-e_{n_L-1})(e_{m_p}-e_{n_L})=e_{m_p}-e_{n_L},$$ it follows that
$(e_{m_p+1}-e_{n_L-1})\sum_{j=L}^{p}d_j= \sum_{j=L}^{p}d_j$. But then $m_p+1>  N$, and hence,
\ba \sum_{j=L}^{p}d_j&= (e_{m_p+1}-e_{n_L-1})\sum_{j=L}^{p}d_j (e_{m_p+1}-e_{n_L-1})\\
&\le (e_{m_p+1}-e_{n_L-1})D (e_{m_p+1}-e_{n_L-1})\\
&\preceq  (e_{m_p+1}-e_{N})D (e_{m_p+1}-e_{N})\\
&\preceq a.
\end{align*}
This proves that $\{d_j\}$ is thin.

\ep

\bT{T:thin}
Let $\A$ be simple, $\sigma$-unital, non-unital, non-elementary C*-algebra with an approximate identity $\{e_n\}$. Then the following are equivalent.
\item [(i)]  $\A_+\ne \overline{K_o(\{e_n\})}$;
\item [(ii)] $\A$ has an order dense sequence;
\item [(iii)] $\A$ has a thin sequence;
\item [(iv)] $\A$ has a thin sequence  $d_j$ such that  $D= \sum_{j=1}^\infty 
d_j$ converges strictly to an element  $D\in K_o(\{e_n\})\setminus \A$.
\eT
\bp
As usual, set $K_o=K_o(\{e_n\})$.\\
(i)$~\Rightarrow~$(ii)  $\A_+\ne \overline{K_o}$ if and only if there is an element $X\in K_o\setminus A$.
Then for every $k$, $(1-e_k)X(1-e_k)\ne 0$, hence there is some integer $m_k>k$ such that $s_k:=(e_{m_k}-e_k)X(e_{m_k}-e_k)\ne 0$. By the defining property of $K_o$, for every $0\ne a \in \A_+$ there is an integer $N$ such that $s_N\preceq a$. \\
(ii)$~\Rightarrow~$(iii) by  Corollary \ref {C: infinit vs thin}\\
(iii)$~\Rightarrow~$(iv) 
Assume that $t_j$ is a thin sequence for $\A_+$. By Lemma \ref {L:smaller elements}, for every $j$ we can find  $0\ne \tilde d_j\in \A_+$ such that $\tilde d_j\preceq t_j$ and $\tilde d_j\le e_{2j}-e_{2j-1}$. Let 
$d_j:= \frac{\tilde d_j}{\|\tilde d_j\|}$ and $D:=\sum_{j=1}^\infty d_j$. The sequence $\{d_j\}$ is mutually orthogonal and thin by  (\ref {e: subthin}) and (\ref {e: subthin2}), and by construction, $D$ is diagonal with respect to $\{e_n\}$.  Then by Lemma \ref {L: diag},  $D\in K_o\setminus \A$.  \\
(iv)$~\Rightarrow~$(i) Obvious.

\ep
Immediate consequences of Theorem \ref {T:thin}, Proposition \ref {P:separable} and Proposition \ref {P: countable dim group}, and Lemma \ref {L:cont scale}  are the following ((i) was obtained in \cite [Lemma 2.4] {LinContScale}). \bC{C:separable} Let $\A$ be a simple, $\sigma$-unital, non-unital, non-elementary C*-algebra with an approximate identity $\{e_n\}$. Then $\A_+\ne \overline {K_o}$ in any of the following cases:
\item [(i)]  $\A$ is separable;
\item [(ii)] the Cuntz semigroup is order separable;
\item [(iii)] $\A$ has the (SP) property and its dimension semigroup $D(\A)$ of Murray-von Neumann equivalence classes of projections is order separable;
\item [(iv)] $\A$ has a continuous scale.
\eC
We will see in Section \ref{S:strict comp} that another case when $\A_+\ne \overline {K_o}$ is when $\A$ has strict comparison of positive elements by traces (see Proposition \ref {P:Icon not A} and Theorem \ref {T: Imin=Icont}).

\section{The minimal ideal}

We proceed now to prove that for every approximate identity $\{e_n\}$, as usual, satisfying (\ref {e: appr id}), $K_o(\{e_n\})$ is weakly invariant and to obtain properties of $\overline{L(K_o(\{e_n\}))}$. In order to do that, we need first to strengthen a result obtained in \cite [Theorem 4.2]{KNZCompPos}. Diagonal series have proven a very valuable tool in working with multiplier algebras, started with \cite {ElliottAnn} and used by \cite {LinThesis},\cite{SZ90},\cite {RordamIdeals} among many other. It is well known that a Weyl-von Neumann decomposition of selfadjoint elements into the sum of a diagonal series plus a remainder in $\A$ of arbitrarily small norm is possible only under additional conditions on $K_1(\A)$ (e.g., if $\A$ has real rank zero, the Weyl-von Neumann theorem holds precisely when $\Ma$ has real rank zero \cite {SZ91}.) However a decomposition into a tridiagonal series plus remainder was obtained and used in \cite {SZ90}.
A  refinement of that construction, but with fewer hypotheses on $\A$, was obtained in \cite {KNZCompPos} where we proved that if $\A$ is $\sigma$-unital, then every positive element  $T\in \Ma_+$ can be decomposed into the sum of a selfadjoint element in $\A$ of arbitrarily small norm and a {\it bidiagonal} series. A bidiagonal series $D:=\sum_{k=1}^\infty d_k$  is a strictly converging series with summands $d_k\in \A_+$ such that $d_kd_{k'}=0$ for $|k-k'|>1$.  In particular, $D=D+e+D_o$, where $D_e:= \sum_{k=1}^\infty d_{2k}$ and $D_o:= \sum_{k=1}^\infty d_{2k-1}$ are diagonal series.

If $T\in K_o(\{e_n\})$, the original proof in \cite {KNZCompPos} can be modified to show that the bidiagonal series can be chosen in $\overline {K_o(\{e_n\})}$.  Also in order to obtain some further enhancements that will be needed later in this paper, and for the readers' convenience, we will present here a self-contained proof. 
\bT {T:bidiag revisited}
Let $\A$ be a simple, $\sigma$-unital, non-unital, non-elementary C*-algebra with approximate identities $\{e_n\}$ and $\{f_m\}$, and  let $X^*X\in K_o(\{e_n\})$ for some $X\in \Ma$. Then for every $\eps>0$, there exist an element $t=t^*\in \A$ with $\|t\|< \eps$, and a bidiagonal series $D:=\sum_{k=1}^\infty d_k$ such that $XX^* =D + t$ and  $D\in K_o(\{f_m\})$.
\eT
\bp 
Without loss of generality, assume that $\|X\|=1$ and assume also that $X^*X\not \in \A$ as the conclusion is trivial when $XX^*\in \A$ (e.g., see (\ref {e: Asubclos K0})). By Theorem  \ref {T:thin}, there exists a thin sequence $t_k$. By the definition of $ K_o(\{e_n\})$ there is an increasing sequence $N_k$,  such that
$$ (e_m-e_{N_k})X^*X(e_m-e_{N_k})\preceq t_{k+1}\quad \forall m>N_k.$$
Since $K_o(\{e_n\})=K_o(\{e_{N_k}\})$ by Lemma \ref {L:Ko for subseq},   to simplify notations assume that
\be {e: prec cn}(e_m-e_{n})X^*X(e_m-e_{n})\preceq t_{n+1} \quad \forall m>n.\ee
Fix $\epsilon>0$ and construct two  sequences of $\{m_k\}$ and $\{n_k\}$ of strictly increasing integers as follows. Set $m_0=n_0=n_{-1}=0$, $n_1=1$ and $e_0=f_0=0.$
Since $\{f_m\}$ is an approximate identity, we can find $m_1>0$ such that 
$$
\|e_{n_1}X^*(1-f_{m_1})\|< \frac{\epsilon^2}{4^{3}}.
$$
Then choose $n_2>n_1=1$ such that
$$\|(1-e_{n_2})X^*f_{m_1}\| < \frac{\epsilon^2}{4^{5}}.$$
By iterating, construct strictly increasing sequences of integers $m_k$ and $n_k$ such that 
\ba  \|e_{n_k}X^*(1-f_{m_k})\|&< \frac{\epsilon^2}{4^{k+2}}\quad\text{for }k\ge 1\\
\|(1-e_{n_{k-1}})X^*f_{m_{k-2} }\| &< \frac{\epsilon^2}{4^{k+2}} \quad\text{for }k\ge 3.
\end{align*}
When  $A,B, C$ are bounded operators, $\|C\|\le 1$, and $0\le A\le B$, then $$\|A^{1/2}C\|^2= \|C^*AC\|\le \|C^*BC\| \le\|BC\|.$$ Using the fact that $\|X\|= 1$ and $\|f_{m_k}\|=1$ for all $k$, we can apply  this inequality to $A:=(e_{n_{k}}-e_{n_{k-1}})$ and 
$$ B:=e_{n_{k}}\quad \text {and}\quad C:=X^*(1-f_{m_{k} })$$
and also to 
$$ B:=1- e_{n_{k-1}}\quad \text {and}\quad C:=X^*f_{m_{k-2} }.$$ 
Thus we obtain
\ba
&\|(e_{n_{k}}-e_{n_{k-1}})^{1/2}X^*(1-f_{m_{k} })\|\le \frac{\epsilon}{2^{k+2}}\quad\text{for }k\ge 1\\
&\|(e_{n_{k}}-e_{n_{k-1}})^{1/2}X^*f_{m_{k-2} }\|\le \frac{\epsilon}{2^{k+2}}\quad\text{for }k\ge 3.\end{align*}
By the triangle inequality,
\ba
\|(e_{n_{k}}-e_{n_{k-1}})^{1/2}&X^*- (e_{n_{k}}-e_{n_{k-1}})^{1/2}X^*\big(f_{m_k}-f_{m_{k-2}}\big)\|\\
&= \|(e_{n_{k}}-e_{n_{k-1}})^{1/2}X^*(1-f_{m_{k} })+ (e_{n_{k}}-e_{n_{k-1}})^{1/2}X^*f_{m_{k-2} }\|\\
&< \frac{\epsilon}{2^{k+1}}.
\end{align*}
From the inequality $\|A^*A-B^*B\|\le (\|A\|+\|B\|)(\|A-B\|)$ and again using the fact that $\|X\|= 1$ and $\|e_{n_{k}}\|=\|f_{m_{k}}\|=1$, we thus have 
\be{e:ineq for Yk}\|X(e_{n_{k}}-e_{n_{k-1}})X^*- \big(f_{m_k}-f_{m_{k-2}}\big)X(e_{n_{k}}-e_{n_{k-1}})X^*\big(f_{m_k}-f_{m_{k-2}}\big)\|\le \frac{\epsilon}{2^{k}}.\ee
Set  
\ba c_k&:= \big(f_{m_k}-f_{m_{k-2}}\big)X(e_{n_{k}}-e_{n_{k-1}})X^*\big(f_{m_k}-f_{m_{k-2}}\big)\\ 
D&:= \sum_{k=1}^\infty c_k.\end{align*}
Since $f_m$ is an approximate identity for $\A$ and the sequence 
$$c_k\le  \|X\|^2 (f_{m_k}-f_{m_{k-2}})^2\le f_{m_k}-f_{m_{k-2}}$$ is uniformly bounded, it is clear that the series converges strictly. Furthermore, 
$$XX^*= \sum_{k=1}^\infty X(e_{n_k}-e_{n_{k-1}} )X^*$$ where the series also convergences strictly. Set 
$$
t:= XX^*-D = \sum_{k=1}^\infty \big(X(e_{n_k}-e_{n_{k-1}} )X^*-c_k\big).
$$ It follows from (\ref {e:ineq for Yk}) that this series converges in norm, hence $t=t^*\in \A$. Moreover 
$$ \|t\| \le \sum_{k=1}^\infty\|X(e_{n_k}-e_{n_{k-1}} )X^*-c_k\| < \epsilon.$$ 
Thus we have the decomposition $XX^*=D+t$.  We need to verify that $D$ is a bidiagonal series and that  $D\in K_o(\{f_m\})$.
We will use now (\ref{e: prec cn}), which is a consequence of  $X^*X\in K_o(\{e_n\})\setminus \A$. For every $k>1$ 
\begin{alignat*}{2}
c_k&\preceq X(e_{n_{k}}-e_{n_{k-1}})X^* & (\text{by Lemma \ref {L:Cuntz subeq} (i)})\notag\\
&\sim (e_{n_{k}}-e_{n_{k-1}})^{1/2}X^*X(e_{n_{k}}-e_{n_{k-1}})^{1/2}\qquad & (  \text{by Lemma \ref {L:Cuntz subeq} (ii)})\notag\\
&\sim (e_{n_{k}}-e_{n_{k-1}})X^*X(e_{n_{k}}-e_{n_{k-1}})& (\text{by Lemma \ref {L:sqrt} (ii)})\notag\\
&\preceq t_{k}&(\text{by  (\ref {e: prec cn})}\label{e:Yk}.
\end{alignat*}
Set $d_k:=c_{2k}+c_{2k-1}$. By Lemma \ref {L:Cuntz subeq} (vi),
\be{e:dk vs tk} d_k\preceq t_{2k}+t_{2k-1}.\ee
Furthermore,
\be{e:dk} d_k\le 2 \|X\|^2 (f_{m_{2k}}-f_{m_{2k-3}})\ee 
whence we see that $D$ is bidiagonal.
In particular, the even and odd sequences
\ba d_{2k}&\le 2 \|X\|^2 (f_{m_{4k}}-f_{m_{4k-3}})\\
d_{2k+1}&\le 2 \|X\|^2 (f_{m_{4k+2}}-f_{m_{4k-1}})
\end{align*}
are both mutually orthogonal, satisfy the intertwining condition  of Lemma \ref {L: diag}, and are thin  by  (\ref{e:dk vs tk}) and (\ref {e: subthin}) since \ba d_{2k}&\preceq  t_{4k}+t_{4k-1}\\
d_{2k+1}&\preceq t_{4k+2}+t_{4k+1}\end{align*} and both sequences $\{ t_{4k}+t_{4k-1}\}$ and $\{ t_{4k+2}+t_{4k+1}\}$ are thin.
But then their sums 
$$D_e:= \sum_{k=1}^\infty d_{2k}\quad\text{and}\quad D_o:= \sum_{k=1}^\infty d_{2k-1}$$
are both in $K_o(\{f_m\})$, and hence, $D=D_e+D_o\in K_o(\{f_m\})$, which concludes the proof.
\ep
\bR{R:bidiag}
If in Theorem \ref {T:bidiag revisited} we start with an element $B\in \Ma_+$ and drop the hypothesis that $B\in K_o(\{e_n\})$, the same proof yields the decomposition $B= D+t$ where $D$ is a bidiagonal series. Furthermore, if $\{f_m\}$ is an approximate identity, then we can choose $D$  to be the sum $D=D_e+D_o$ of two diagonal series with respect to $\{f_m\}$.  In fact to obtain this result we only need to require that $\A$  is $\sigma$-unital (see \cite[Theorem 4.2] {KNZCompPos}).
\eR
\bC {C: Imin} 

Let $\A$ be simple, $\sigma$-unital, non-unital, non-elementary and let $e_n,\, f_m$ be two approximate identities for $\A$. Then
\item [(i)]  $\overline{K_o(\{e_n\})}=\overline{K_o(\{f_m\})}$
\item [(ii)] $K_o(\{e_n\})$ is weakly invariant, hence $\overline {K_o(\{e_n\})}$   is hereditary and strongly invariant.
\item [(iii)] $\overline{ L(K_o(\{e_n\}))}$ is a two-sided ideal and $$\overline{ L(K_o(\{e_n\}))}=L(\,\overline {K_o(\{e_n\}}\,)=\text{span}\{ \overline{K_o(\{e_n\})}\}.$$
\eC
\bp
\item [(i)] 
If $T\in K_o(\{e_n\})$, then by applying Theorem  \ref {T:bidiag revisited} to $X:=T^{\frac{1}{2}}$ we see that $T\in \overline{K_o(\{f_m\})}$, that is,  $K_o(\{e_n\})\subset \overline{K_o(\{f_m\})}$. Thus $\overline{K_o(\{e_n\})}\subset \overline{K_o(\{f_m\})}$. By reversing the role of the approximate identities we obtain equality.
\item [(ii)]  If $X\in  K_o(\{e_n\})$ and $A\in \Ma$ then $$(X^{1/2}A^*)(X^{1/2}A^*)^*
=X^{1/2}A^*AX^{1/2}\le \|A\|^2X\in K_o(\{e_n\}), $$ hence by Theorem  \ref {T:bidiag revisited}, $$AXA^*=(X^{1/2}A^*)^*(X^{1/2}A^*)\in \overline{K_o(\{e_n\})}.$$
Thus $K_o(\{e_n\})$ is weakly invariant, and hence, by Lemma \ref {L:new stuff cone}(ii) and (iv), we obtain that $\overline {K_o(\{e_n\})}$   is hereditary and strongly invariant.
\item [(iii)] Follows immediately from Corollary \ref {C:closure ideal} and Theorem \ref {T cones}.
\ep
The independence of $\overline{ L(K_o(\{e_n\}))}$ on the approximate identity was obtained  in \cite[Remark 2.9]{LinContScale}. From now on, we will denote

\be{e:def Imin}\Imin:= \overline{ L(K_o(\{e_n\}))}.\ee

The following result sheds some additional light on the relation between $\Imin$ and $ L(K_o(\{e_n\}))$.

\bP{P: Imin and Io}
Let $\A$ be simple, $\sigma$-unital, non-unital, non-elementary and let $\{e_n\}$ be an approximate identity for $\A$. Then $\Imin= \A+ L(K_o(\{e_n\})).$
\eP
\bp
The inclusion $\A + L(K_o(\{e_n\}))\subset \Imin$ is obvious, and to prove equality it is enough to verify that if $D\in (\Imin)_+= \overline{K_o(\{e_n\})} $, then $D\in  \A+ L(K_o(\{e_n\}))$. 

Without loss of generality, $\|D\|\le 1$ and by Remark \ref {R:bidiag} we can assume that $D$ is diagonal with respect to $\{e_n\}$. By further decomposing if necessary $D=\sum_{j=1}^\infty d_j$ into a sum of at most three diagonal series, we can assume that there is a sequence $m_k$ such that $(e_{m_k}-e_{m_{k-1}})d_k= d_k$ for all $k$. To simplify notations, assume that 
\be{e: isolated} (e_k-e_{k-1})d_k=d_k\quad \forall k\ee (setting $e_0=0$). 
By Theorem \ref {T:thin}, $\A$ has a thin sequence $\{t_j\}$. For every $k$ find $b_k\in K_o(\{e_n\})$ such that $\|D-b_k\|< \frac{1}{k}$ and an integer $n_k$ such that $$(e_m-e_{n_k})b_k(e_m-e_{n_k})\preceq t_k\quad \forall m >n_k.$$ Since
$$\| (e_m-e_{n_k})D(e_m-e_{n_k})-(e_m-e_{n_k})b_k(e_m-e_{n_k}) \| \le \|D-b_k\|<  \frac{1}{k}$$
it follows from Lemma \ref {L:Cuntz subeq} (iv) that for all $m>n_k$,
$$\big ((e_m-e_{n_k})D(e_m-e_{n_k})-\frac{1}{k}\big)_+\preceq (e_m-e_{n_k})b_k(e_m-e_{n_k})\preceq t_k.$$
By (\ref {e: isolated}),
$$\big ((e_m-e_{n_k})D(e_m-e_{n_k})-\frac{1}{k}\big)_+= \Big( \sum_{j=n_k}^m d_j-\frac{1}{k}\Big)_+ = \sum_{j=n_k}^m (d_j-\frac{1}{k})_+.$$
Set $\delta_j: =\frac{1}{k} $ for $n_{k}\le j< n_{k+1}$. Thus for all $k\in \mathbb N$,
$\sum_{j=n_k}^{n_{k+1}-1} (d_j-\delta_j)_+\preceq t_k$. Then for every $0\ne a\in \A_+$ there is an $K\in \mathbb N$ such that $\sum_{j=K}^k t_j \preceq a$ for all $k>K$. For all $m >n_K$, choose $n_H\ge m$. Then 
\ba(e_m- e_{n_K})&\Big(\sum_{j=1}^\infty (d_j-\delta_j)_+\Big)(e_m- e_{n_K})= \sum_{j=n_K}^m (d_j-\delta_k)_+\\
&\le\sum_{j=n_K}^{n_H} (d_j-\delta_k)_+\le  \sum_{k=K}^{H-1} \sum_{j=n_k}^{n_{k+1}} (d_j-\delta_j)_+\preceq \sum_{k=K}^{H-1}t_k\preceq a\end{align*}
which proves that $$\sum_{j=1}^\infty (d_j-\delta_j)_+ \in K_o(\{e_n\})\subset L( K_o(\{e_n\})).$$ Finally, 
$$D- \sum_{j=1}^\infty (d_j-\delta_j)_+= \sum_{j=1}^\infty (d_j-(d_j-\delta_j)_+)\in \A_+$$
since $0\le d_j-(d_j-\delta_j)_+\le \delta_j(e_{j+1}-e_j)$.
\ep

We proceed now to justify the notation $\Imin$. The natural  ``minimal ideal" is  the intersection $\Io$ of all ideals (not necessarily proper) properly containing $\A$, in symbols

\be {e:Ioo}\Io:=\bigcap\{ \J\lhd \Ma, \A\subsetneq\J\}.\ee
Obviously $\A\subset \Io$, but we do not know whether $\A\ne \Io$ holds in general.  However, we will prove now that $\Imin=\Io$ (see Theorem \ref {T: minimal}). A key tool in that proof, and used also throughout this paper, is the following result obtained in \cite{{KNZCompPos}}.
\bP{P:diag Cuntz} \cite[Proposition 3.2]{KNZCompPos} Let $\B$ be a non-unital  C*-algebra and let $A=\sum_{n=1}^\infty A_n$, $B=\sum_{n=1}^\infty B_n$ where $A_n, B_n\in \Mul(\B)_+$, $A_nA_m=0$, $B_nB_m=0$ for $n\ne m$ and the two series converge in the strict topology, and $A_n\preceq (B_n-\delta)_+$ for some $\delta >0$ and for all $n$. Then for every $\epsilon>0$ and $0< \delta'< \delta$ there is an $X\in \Mul(\B)$ such that $(A-\epsilon)_+=X(B-\delta')_+X^*$, and hence, $A\preceq (B-\delta')_+\le B$. 
\eP
If the sum of a positive diagonal series in $\Mb$ is subequivalent to another strictly converging series in $\Mb$ (not necessarily diagonal) then we can deduce the following relations between the summands.

\bP{P: from subeq} Let $\B$ be a non-unital C*-algebra, $A= \sum_{k=1}^\infty a_k$,  $B~=~\sum_{k=1}^\infty b_k$, be two  strictly converging series with $a_k,\,b_k\in \B_+$ and elements $a_k$ mutually orthogonal. If $A\preceq B$, then for every $\delta >0$ and $M\in \mathbb N$ there is an $N\in \mathbb N$ such that for every $n\ge N$ there is an $m\ge M$ such that
$$\sum_{k=N}^n(a_k-\delta)_+\preceq \sum_{k=M}^mb_k.$$
\eP
\bp
By Lemma \ref {L:Cuntz subeq} (v), there is an $X\in \Mb$ such that $(A-\frac{\delta}{6})_+=XBX^*$, and hence, by Lemma \ref {L: approx} there is a $Y \in \Mb$ such that \be{e:1/6} \|(A-\frac{\delta}{6})_+- \big((A-\frac{\delta}{6})_+\big)^{1/2}YXBX^*Y^*\big((A-\frac{\delta}{6})_+\big)^{1/2}\|< \frac{\delta}{6}\ee and $\|YXBX^*Y^*\|\le 1$.
Because of the mutual orthogonality of  $a_k$, and hence, of $(a_k-
\frac{\delta}{6})_+$, we have for every $n$
\be{e: trivial ineq}(A-\frac{\delta}{6})_+= \sum_{k=1}^\infty (a_k-
\frac{\delta}{6})_+\ge \sum_{k=n}^\infty (a_k-
\frac{\delta}{6})_+.\ee
If  $a, b,c$ are positive elements in a C*-algebra $\C$ with $a\le b$ and $\|c\|\le 1$,  then 
\begin{align}\label {e:minilemma}
\|a-a^{1/2 }ca^{1/2 }\|&= \|a^{1/2 }(1-c)a^{1/2 }\|
=\|(1-c)^{1/2}a(1-c)^{1/2}\|\\
&\le \|(1-c)^{1/2}b(1-c)^{1/2}\|= \|b-b^{1/2 }cb^{1/2 }\|.\notag
\end{align}
Thus from (\ref{e:1/6}), (\ref {e: trivial ineq}), and (\ref {e:minilemma}) we have for all $n$
\be{e: step1}
\| \sum_{k=n}^\infty(a_k-\frac{\delta}{6})_+- \big( \sum_{k=n}^\infty(a_k-\frac{\delta}{6})_+\big)^{1/2}YXBX^*Y^*\Big(\sum_{k=n}^\infty(a_k-\frac{\delta}{6})_+\Big)^{1/2}\|< \frac{\delta}{6}.
\ee
Since $YX\sum_{k=1}^{M-1}b_kX^*Y^*\in \B$ and  the sequence $\sum_{k=n}^\infty (a_k-\frac{\delta}{6})_+ \to 0$ strictly, we can find an integer $N$ such that
\be{e:2/6}\| \Big(\sum_{k=N}^\infty(a_k-\frac{\delta}{6})_+\Big)^{1/2}YX\sum_{k=1}^{M-1}b_kX^*Y^*\Big(\sum_{k=N}^\infty(a_k-\frac{\delta}{6})_+\Big)^{1/2}\|< \frac{\delta}{6}.\ee
As a consequence of   (\ref {e: step1}) and (\ref {e:2/6}) we thus obtain
\be{e:3/6}\| \sum_{k=N}^\infty(a_k-\frac{\delta}{6})_+- \big( \sum_{k=N}^\infty(a_k-\frac{\delta}{6})_+\big)^{1/2}YX\sum_{k=M}^{\infty}b_kX^*Y^*\Big(\sum_{k=N}^\infty(a_k-\frac{\delta}{6})_+\Big)^{1/2}\| <  \frac{2\delta}{6},\ee
and hence,
$$\sum_{k=N}^\infty(a_k-\frac{3\delta}{6})_+\preceq \big( \sum_{k=N}^\infty(a_k-\frac{\delta}{6})_+\big)^{1/2}YX\sum_{k=M}^{\infty}b_kX^*Y^*\Big(\sum_{k=N}^\infty(a_k-\frac{\delta}{6})_+\Big)^{1/2}\preceq \sum_{k=M}^{\infty}b_k.$$
A fortiori, for every $n\ge N$, we have $\sum_{k=N}^n(a_k-\frac{3\delta}{6})_+\preceq
\sum_{k=M}^{\infty}b_k. $
Then again by  Lemma \ref {L:Cuntz subeq} (v), there  is a $Z\in \Mb$ such that 
$$\sum_{k=N}^n(a_k-\frac{4\delta}{6})_+ = Z\sum_{k=M}^{\infty}b_kZ^*.$$
Choose $e\in\B$ such that 
$\|Z\sum_{k=M}^{\infty}b_kZ^*- eZ\sum_{k=M}^{\infty}b_kZ^*e\|< \frac{\delta}{6}.$
and then choose $m\ge M$ such that 
$ \|eZ\sum_{k=m+1}^{\infty}b_kZ^*e\|< \frac{\delta}{6}.$
Then
$$\|\sum_{k=N}^n(a_k-\frac{4\delta}{6})_+- eZ\sum_{k=M}^{m}b_kZ^*e\|< \frac{2\delta}{6},$$
and hence
$$\sum_{k=N}^n(a_k-\delta)_+\preceq  eZ\sum_{k=M}^{m}b_kZ^*e\preceq \sum_{k=M}^{m}b_k.$$
\ep

The  inclusion $\Imin \subset \mathcal I_{o}$ in the following theorem has been obtained in \cite [Theorem 2.8] {LinContScale}, but for completeness's sake we include its proof.

\bT{T: minimal}
Let $\A$ be a simple, $\sigma$-unital, non-unital, non-elementary C*-algebra. Then   $\Imin=\Io$ 
\eT 
\bp
To prove that $\Imin \subset \Io$, it is enough to show that given an approximate identity $\{e_n\}$, an element $D\in K_o(\{e_n\})$ and an element $C\in \Ma_+\setminus \A$, then $D\in I(C)$. 
 By Theorem \ref {T:bidiag revisited} and Remark \ref {R:bidiag}, $C= C_e+C_0+ t$ for some $t=t^*\in \A$ and two positive diagonal series $C_e$ and $C_o$ (with respect to $\{e_n\}$), at least one of which, say $C_e$,   does not belong to $\A$. Then, $I(C_e)\subset I(C_e+C_o)=I(C)$, thus it is enough to prove that $D\in I(C_e)$. To simplify notations, assume that $C= \sum_{k=1}^\infty c_k$  itself is diagonal with respect to $\{e_n\}$.
 By Theorem \ref {T:bidiag revisited} and Remark \ref {R:bidiag}, we can also assume  that the series $D= \sum_{k=1}^\infty d_k$ is diagonal with respect to $\{e_n\}$.
Since $\underset{\delta\to 0}{\lim}(C-\delta)_+=C\not \in \A$, there is some $\delta >0$ such that $(C-\delta)_+ \not \in \A$.  Since $(C-\delta)_+= \sum_{k=1}^\infty (c_k-\delta)_+$, we can assume without loss of generality that $(c_k-\delta)_+\ne 0$ for every $k$.

By Lemma \ref {L: diag}, the sequence $\{d_j\}$ is thin, hence  for every $k$ there is an integer $n_k$ such that $$\sum_{j=n_k+1}^m d_j\preceq (c_k-\delta)_+ \quad \forall \, m\ge n_k, ~ k\in \mathbb N.$$
Choose the sequence $n_k$  so to be strictly increasing.
Then in particular
$$
\sum_{j=n_k+1}^{n_{k+1}} d_j\preceq (c_k-\delta)_+ \quad \forall \, k\in \mathbb N,$$
and hence, by Proposition \ref {P:diag Cuntz}, 
$$
\sum_{j=n_1+1}^\infty  d_j\preceq (C-\frac{\delta}{2})_+\le C.$$
Thus  $ D\in I(C)$,  which shows that $\Imin \subset \Io$.

 Now to prove that $\Io=\Imin$, we need to consider only the case  that  $\A\ne \Io$. We will prove  that then   $\Io$  contains a thin sequence, which  by Theorem \ref {T:thin} implies that $\A\ne \Imin$ and hence that $\Io\subset \Imin$.  Equality then holds by the first part of the proof.

Choose $D\in (\Io)_+\setminus \A$ and by invoking  Theorem \ref {T:bidiag revisited} and Remark \ref {R:bidiag} as in the first part of the proof, assume that $D:=\sum_{k=1}^\infty d_k$ is diagonal with respect to $\{e_n\}$. 
 Let $\delta>0$ be such that $(D-\delta)_+\not \in \A$. We claim that the sequence $\{(d_k-\delta)_+\}$ is thin. Since $\sum_{k=1}^\infty (d_k-\delta)_+= (D-\delta)_+\not \in \A$, we can assume without loss of generality that $(d_k-\delta)_+\ne 0$ for all $k$. Let $0\ne a\in \A_+$. 
By Lemma \ref {L:ncopies} applied to the stationary sequence $s_i=a$, there is a sequence  of mutually orthogonal elements $0\ne t_i\in \A_+$ such that $n\sum_{i=n}^mt_i\preceq a $ for every pair of integers $m\ge n$. 
By Lemma \ref {L:smaller elements} there are elements $0\ne a'_i\le e_{2i}-e_{2i-1}$ and $a'_i\preceq t_i$ for every $i$. Let $a_i:=\frac{a_i'}{\|a_i'\|}$. Then the series converges strictly to an element $A:= \sum_{i=1}^\infty a_i \in\Ma\setminus  \A$ because $a_i\le \|a_i'\|(e_{2i}-e_{2i-1})$ and  $\|a_i\|=1$ for every $i$. Furthermore, for every $m\ge M\in \mathbb N$ we have
\be{e: A is small} M\sum_{i=M}^ma_i\sim M\sum_{i=M}^m a'_i\preceq M\sum_{i=M}^m t_i\preceq a.\ee
Since $\A\subsetneq I(A)$, it follows that $\Io\subset I(A)$, and hence, there is some $M$ such that $(D-\frac{\delta}{2})_+\preceq MA$. By Proposition \ref {P: from subeq},   there is some $N$ such that for every $n\ge N$ there is $m\ge M$ for which
$$\sum_{k=N}^n(d_k-\delta)_+\preceq \sum_{i=M}^m Ma_i\sim M\sum_{i=M}^m a_i\preceq a.$$
This proves that the sequence $(d_k-\delta)_+$ is thin and thus concludes the proof.

\ep

In \cite {LinSimple}, Lin proved that if $\Ma/\A$ is simple, then it is purely infinite. Thus if $\Imin= \Ma$, and hence $\Ma/\A$ is simple, then $\Imin/\A$ is purely infinite. We can relax the condition $\Imin= \Ma$.

\bT{T:purely infinite} Let $\A$ be a simple, $\sigma$-unital, non-unital, non-elementary C*-algebra and assume that $\Imin\ne \A$. Then $\Imin/\A$ is purely infinite simple.
\eT
\bp By  Theorem \ref {T: minimal}, it is trivial to see that $\Imin/\A$ is simple.
Denote by $\pi:~ \Imin\to \Imin/\A$ the canonical quotient map. Choose a positive element $T\in \Imin\setminus \A$. Given an approximate identity $\{e_n\}$, by Theorem \ref{T:bidiag revisited} and Remark \ref {R:bidiag} we can find a series $D:=\sum_{k=1}^\infty d_k$ diagonal with respect to $\{e_n\}$ and with $0\ne \pi(D)\le \pi(T)$. Choose $\delta>0$ such that $(D-\delta)_+\not\in \A$. By the diagonality of $D$,  $(D-\delta)_+= \sum_{k=1}^\infty (d_k-\delta)_+$ and  assume that $(d_k-\delta)_+\ne 0$ for every $k$. Apply Lemma \ref {L:ncopies} to the sequence
$\{(d_k-\delta)_+\}$ to find a mutually orthogonal sequence $0\ne c''_k\in \A_+$ such that  $nc''_k\preceq (d_k-\delta)_+$ for every $n\in \mathbb N$ and $k\ge n$, where $nc''_k$ denotes as before the n-fold direct sum of $c''_k$ with itself. Choose $0\ne c'_k\le e_{2k}-e_{2k-1}$ with $c'_k\preceq c''_k$ for every $k$. Define $c_k:= \frac{c'_k}{\|c'_k\|}$ and $C:= \sum_{k=1}^\infty c_k$. Then the series converge strictly and $C\not \in \A$. Moreover, $$nc_k\preceq nc''_k\preceq (d_k-\delta)_+\quad  \forall\, k\ge n.$$ By Proposition \ref {P:diag Cuntz}, $$n\sum_{k=n}^\infty c_k\preceq \sum_{k=n}^\infty d_k.$$
But then $$n\pi(C) = n \pi\big(\sum_{k=n}^\infty c_k\big)  \preceq \pi\big(\sum_{k=n}^\infty d_k\big)= \pi(D)\le \pi(T)\quad\forall\, n\in \mathbb N.$$
In particular, $\pi(C)\preceq \pi(T)$, that is, $C\in (\Imin)_+\setminus \A$. 

On the other hand,  $\Imin/\A= \Io/\A$ by Theorem  \ref {T: minimal}, and hence, it is simple.  Thus for every $\epsilon>0$  there is an $m$ such that $\pi((T-\epsilon)_+)\preceq m\pi(C)$, and hence, 
$$\pi((T-\epsilon)_+)\oplus \pi((T-\epsilon)_+)\preceq 2m\pi(C)\preceq \pi(T).$$
Since $\epsilon $ is arbitrary, it follows that $\pi(T)\oplus\pi(T)\preceq \pi(T)$ which proves that $\Imin/\A$ is purely infinite.

\ep

\section{The minimal ideal  when $\A$ has strict comparison}\label{S:strict comp}

\bD{D:Icon}
Let $\A$ be a simple, $\sigma$-unital, non-unital C*-algebra with nonempty tracial simplex $\TA$. Set:  
\item [(i)] $K_c:=\{X\in \Ma_+\mid \hat X\in \Aff\}$.
\item [(ii)] $\Icon:= \overline {L(K_c)}$.
\eD
\bP {P: Kc is hered cone}
Let $\A$ be a simple, $\sigma$-unital, non-unital C*-algebra with nonempty tracial simplex $\TA$. Then
\item [(i)] $K_c$ is a hereditary strongly invariant cone; $L(K_c)$ is a two-sided selfadjoint ideal and hence so is $\overline {L(K_c)}= L(\overline{K_c})$.
\item [(ii)]  \begin{align*}(\Icon)_+&=\overline{K_c}=\{X\in \Ma_+\mid \hat X\in \Aff \, \}\overline{\phantom{m}}\\
&=\{X\in \Ma_+\mid \widehat {(X-\delta)_+}\in \Aff ~\forall ~\delta>0 \}\end{align*}
\item [(iii)] For a projection  $P\in \Ma$,  $P\in \Icon$ if and only if $\hat P$ is continuous.
\item [(iv)] $\Icon= \overline{\spn K_c}$.
\eP
\bp \item [(i)]  
Since the map $\Ma_+\ni X\to \hat X\in \LAff_+$ satisfies the conditions 
$\widehat{X+Y}= \hat X+ \hat Y$  and $\widehat{tX}= t\hat X$ for $X,Y\in \Ma_+$ and $t\in \mathbb R_+$,   it is clear that $K_c$ is a cone. Moreover, if  $0\le X\le Y\in K_c$, then  
$$\hat X+ \widehat {Y-X}=\hat Y.$$
Since $\hat Y$ is affine and continuous  and both $\hat X$ and $ \widehat {Y-X}$ are affine, lower semicontinuous, and non-negative,  it is immediate to verify that both must be continuous. Thus $X\in K_c$,  and hence, $K_c$ is hereditary. Since $\widehat {X^*X}= \widehat {XX^*}$ for all $X\in \Ma$, $K_c$ is strongly invariant.
Therefore, the rest of the conclusions in (i) follows from (\ref {e:str invar}), Lemma \ref {L:new stuff cone}, (\ref{e:invar}), and Corollary \ref  {C:closure ideal}.
\item [(ii)] By Corollary \ref{C:closure ideal}  and Theorem \ref {T cones} (i) and (ii)  we have that $$\big(\overline{L(K_c)}\big)_+= L(\overline{K_c})_+= \overline{K_c}$$ 
which is the first equality in (ii). The second equality is given by Lemma \ref {L:new stuff cone} (iii) 
\item [(iii)] Since $(P-\delta)_+= \begin{cases}(1-\delta)P&0\le \delta< 1\\0&\delta\ge1\end{cases}$, we have by (ii) that $P\in(\Icon)_+$ if and only if $\hat P\in \Aff$.
\item [(iv)] Since by (i) and Theorem \ref {T cones}, $ \Icon= L(\overline{K_c}) = \spn \overline{K_c}$ is closed, it is immediate to see that  $ \spn \overline{K_c}=  \overline{\spn K_c}.$ 
\ep
Notice that if $\A=\K$, then $K_c$ consist of the positive cone of the trace class operators, and hence, $\Icon=\K$.

It is immediate to verify that $\A \subset \Icon$. Indeed  $(a-\delta)_+\in \Ped$ for every $\delta >0$ and  $a\in \A_+$, hence $\widehat{ (a-\delta)_+}$ is continuous, that is, $(a-\delta)_+ \in K_c$. Thus $a\in \overline{K_c}\subset  \Icon$.
To further relate $\Icon$  to $\A$ and to $\Imin$ we need first the following  lemma.

\bL{L: elems with small dim} 
Let $\A$ be a simple, non-elementary  C*-algebra with $\TA\ne \emptyset$. Then for every $\eps>0$ there is an element $0\ne c\in \A_+$ such that $d_\tau(c)< \eps$ for all $\tau\in \TA$. Furthermore, the element $c$ can be chosen in $\Ped$.
\eL
\bp 
Let $0\ne f\in \Ped_+$ and recall that $\hat f\in \Aff_+$.  Choose $\delta>0$ such that $(f-\delta)_+\ne 0$, and an integer $M\ge \frac{\|\hat f\|}{\eps\delta}$. 
By Lemma \ref {L:direct sums} we can find nonzero positive mutually orthogonal  
elements $a_j$ such that $\sum_{j=1}^M a_j\preceq (f-\delta)_+$. By Lemma \ref {L:smaller elements} choose a nonzero positive element $c\preceq a_j$ for $1\le j\le M$. By Lemma \ref {L:Cuntz subeq} (vi) it follows that 
$$M c\preceq  \sum_{j=1}^M a_j\preceq (f-\delta)_+.$$ Thus for every $\tau\in \TA$
\begin{alignat*}{2}Md_\tau(c)&= d_\tau( M c )&(\text{by (\ref {e: dim sum})})\\
&\le d_\tau((f-\delta)_+) \qquad&(\text{by (\ref {e: suub->ineq})})\\
&\le \frac{1}{\delta} \tau(f) \qquad&(\text{by (\ref {e:ineq 2})})\\
&\le  \frac{1}{\delta} \|\hat f\|.
\end{alignat*}
Thus $d_\tau(c)< \eps$. Finally, $(c-\delta)_+\in \Ped$ for every $\delta>0$.  Choose $\delta>0$ such that $(c-\delta)_+\ne 0$. Then $d_\tau((c-\delta)_+)\le d_\tau(c)< \eps$ for all $\tau\in \TA$.
\ep

\bP{P:Icon not A}
Let $\A$ be a simple, $\sigma$-unital, non-unital, non-elementary C*-algebra with nonempty tracial simplex $\TA$. Then $\A\subsetneq \Icon$.
\eP
\bp
 By Lemma \ref {L: elems with small dim}, there is an infinite sequence of elements $0\ne \tilde a_k\in \A_+$ such that $d_\tau(\tilde a_k)\le \frac{1}{2^k}$ for all $k$ and all $\tau\in\TA$. By Lemma \ref {L:smaller elements} we can find $0\ne a'_k\le e_{3k}-e_{3k-1}$ with $a'_k\preceq \tilde a_k$ for all $k$. Let $a_k:=\frac{a'_k}{\|a'_k\|}$. Then $$\tau(a_k)\le  d_\tau(a_k)\le   d_\tau(a'_k)\le d_\tau(\tilde a_k)\le  \frac{1}{2^k}~\forall \,k\text{ and }\forall \,\tau\in\TA.$$ Furthermore, $a_k\le \|a'_k\|(e_{3k+1}-e_{3k-2})\in \Ped$, hence $\widehat {a_k}\in \Aff_+$. Let $A:= \sum_{k=1}^\infty a_k$. Then the series converges strictly  and since it is diagonal ($a_ka_{k'}=0$ for $k\ne k'$) and does not converge in norm, $A\not \in \A$. On the other hand, $\hat A=  \sum_{k=1}^\infty \widehat {a_k}$ is continuous since the series is uniformly convergent. Thus $A\in( \Icon)_+$.
\ep

\bP{P:Imin}
Let $\A$ be a simple, $\sigma$-unital, non-unital,    non-elementary C*-algebra with nonempty $\TA$. Then $K_o(\{e_n\})\subset K_c $ for every   approximate identity $\{e_n\}$. Consequently, $\Imin\subset \Icon$.

\eP
\bp
Let $0\ne X\in K_o(\{e_n\})$ and $\epsilon >0$. By Lemma \ref {L: elems with small dim} we can find an element  $0\ne c\in \A_+$ such that $d_\tau(c)< \frac{\epsilon}{\|X\|}$ for every $\tau\in \TA$.
By the definition of  $K_o(\{e_n\})$ there is an $N\in \mathbb N$ such that $$(e_n-e_m)X(e_n-e_m)\preceq c\quad \forall\,n>m \ge N.$$ Now  $\widehat{X^{1/2}(e_n-e_m) X^{1/2}}\in \Aff_+$ because $X^{1/2}(e_n-e_m) X^{1/2}\in \Ped$. Moreover, 
\begin{align*} 
\widehat{X^{1/2}(e_n-e_m) X^{1/2}}(\tau)&= \widehat{(e_n-e_m)^{1/2} X(e_n-e_m)^{1/2}}(\tau)\\
&\le \|X\|d_\tau\big((e_n-e_m)^{1/2} X(e_n-e_m)^{1/2}  \big)\\
&=  \|X\|d_\tau\big((e_n-e_m) X(e_n-e_m)  \big)\\
&\le  \|X\|d_\tau(c)< \epsilon
\end{align*}
Thus the series $\hat X= \sum_{n=1}^\infty \widehat{X^{1/2}(e_n-e_{n-1}) X^{1/2}}$ converges uniformly and  hence $X\in K_c$. 
This proves that $K_o(\{e_n\})\subset K_c$, and hence, $\Imin\subset \Icon$.

\ep 

In general,  $\Imin$ may fail   to coincide with $\Icon$ as we will see in section \ref{S:Villadsen}.

\bT{T: Imin=Icont}
Let $\A$ be a simple, $\sigma$-unital, non-unital, non-elementary C*-algebra with strict comparison of positive elements by traces. Then  $\Imin= \Icon$.
\eT
\bp
By Proposition \ref {P:Imin} we need to prove that $(\Icon)_+\subset (\Imin)_+$. As in the proof of Theorem \ref{T: minimal}, it is enough to verify that if $\{e_n\}$ is an approximate identity for $\A$, $D=\sum_{k=1}^\infty d_k$ is diagonal with respect to $\{e_n\}$, and $D\in\Icon$, then $D\in \Imin$. Let $\delta>0$, and by dropping if necessary the zero summands in the series $(D-\delta)_+= \sum_{k=1}^\infty(d_k-\delta)_+$, assume that $(d_k-\delta)_+\ne 0$ for all $k$. We claim that the sequence $\{(d_k-\delta)_+\}$ is thin. 

Let $0\ne a \in \A_+$. Recall that the function $d_\tau(a)$ is lower semicontinuous, and hence, $\min_{\tau\in \TA}d_\tau(a)>0$.  By Proposition \ref {P: Kc is hered cone}, $\widehat{(D-\delta)_+}\in \Aff$ and since $(d_k-\delta)_+\in \Ped$ for all $k$, also $\widehat{(d_k-\delta)_+}\in \Aff$. Since
$$\tau((D-\delta)_+)= \sum_{k=1}^\infty \tau((d_k-\delta)_+),$$ by Dini's Theorem the series converges uniformly on $\TA$ for every $\delta >0$.
In particular, there is an $N$ such that if $j\ge i\ge N$ and $\tau\in \TA$, then
\be{e:Cauchy}\sum_{k=i}^j\tau((d_k-\frac{\delta}{2})_+)< \frac{\delta}{2}\min_{\tau\in \TA}d_\tau(a).\ee
By (\ref {e:ineq 2}), $d_\tau( (d_k-\delta)_+) \le \frac{2}{\delta}\tau((d_k-\frac{\delta}{2})_+)$, and hence, by (\ref {e: dim sum}), 
$$d_\tau\big( \sum_{k=i}^j (d_k-\delta)_+\big)= \sum_{k=i}^j d_\tau((d_k-\delta)_+)\le  \sum_{k=i}^j \frac{2}{\delta}\tau((d_k-\frac{\delta}{2})_+)< \min_{\tau\in \TA}d_\tau(a).$$
By the hypothesis of strict comparison of positive elements by traces, we thus have that
$\sum_{k=i}^j (d_k-\delta)_+\preceq a$, which proves that the sequence $\{(d_k-\delta)_+\}$ is thin. But then $(D-\delta)_+\in K_o(\{e_n\})$ by Lemma \ref {L: diag}. Since $\delta$ is arbitrary, it follows that  $D\in \Imin = \overline{K_o(\{e_n\})} $ which concludes the proof.

\ep

As a consequence of this theorem, any counterexample for $\Imin \ne \A$, could only be found among non-separable C*-algebras with no strict comparison of positive elements. Among such algebras is the C*-algebra $\A$ introduced by Rordam to provide an example of a simple unital  C*-algebra  with both infinite projections and nonzero finite projections (\cite [Theorem 5.6] {RordamEx}). Recall that $\A$ is the C*-inductive limit
$\A = \lim_{n \rightarrow \infty} \Mul(\B_n),$
where all $\B_n$ are separable C*-algebras. So, while the algebras  $\Mul(\B_n)$ are not separable and hence neither is $\A$, the order dense sequences for $\B_n$ are order dense also for $\Mul(\B_n)$ and therefore their union is order dense for $\A$. As a consequence $\Imin \ne \A$.

\section{Strict comparison in the minimal ideal.}
In \cite [Theorem 6.6] {KNZCompPos} we proved that if $\A$ is a $\sigma$-unital simple C*-algebra with strict comparison of positive elements by traces and with quasicontinuous scale (e.g.,  with finite extremal boundary), then  strict comparison of positive element by traces (see Definition \ref {D:str comp}) holds also in $\Ma$.  In this section we will show that if we restrict our attention to comparison between elements in $\Icon$, then strict comparison holds {\it without}  requiring the scale to be quasicontinuous. 

For the first step we list here a slightly modified version  of   \cite [Lemma 6.2] {KNZCompPos}. 

\bL{L: abs cont in Imin}
Let $\A$ be a simple, $\sigma$-unital, non-unital C*-algebra with nonempty tracial simplex $\TA$ and let $A\in (\Icon)_+$, $B\in \Ma_+$, and assume that $d_\tau(A)< d_\tau(B)$ for every $\tau\in \TA$ for which $d_\tau(B)< \infty$. Then for every $\epsilon>0$ there is a $\delta>0$ and $\alpha >0$ such that $d_\tau\left(\left(A-\epsilon\right)_+\right) +\alpha \le d_\tau\left(\left(B-\delta\right)_+\right)$ for every $\tau\in \TA$
\eL
The proof being essentially the same, we refer the reader to   \cite [Lemma 6.2] {KNZCompPos}. The only difference is  that here we need to replace the condition used in \cite [Lemma 6.2] {KNZCompPos} that $\hat A\mid_K$ is continuous for some closed subset $K$ of $\TA$, with the condition that $\widehat {(A-\frac{\eps}{2})_+}$ is continuous on the whole of $\TA$, which follows from Proposition \ref {P: Kc is hered cone}.

The next lemma extends the results of \cite [Lemma 6.4] {KNZCompPos}. 
\bL{L: bidiagonal}
Let $\A$ be a simple, $\sigma$-unital, non-unital  C*-algebra with nonempty tracial simplex $\TA$ and let $B=\sum_{k=1}^\infty b_k$ be a strictly converging  series with $b_k\in \A_+$ for all $k$ and $b_nb_m=0$ for $|n-m|\ge 2$. Assume that $B\in (\Icon)_+$ and that $\delta>0$. Then 
\item [(i)]
$d_\tau\left(\left( \sum_n^\infty b_k-\delta\right)_+\right)\downarrow 0$ uniformly on $\TA$.
\item [(ii)] For every $\epsilon >0$ and $0< \delta'< \delta$ there is an $n$ such that for all $\tau\in \TA$
$$d_\tau\left(\left( \sum_{k=1}^n b_k-\delta'\right)_+\right)> d_\tau\left(\left( \sum_{k=1}^\infty b_k-\delta\right)_+\right)-\epsilon.$$
\eL
\bp
\item [(i)] 
The sequence $d_\tau\left(\left( \sum_n^\infty b_k-\delta\right)_+\right)$ is monotone decreasing by Lemma \ref {L:Cuntz subeq} (viii) and (\ref {e: suub->ineq}). Moreover, by Lemma \ref {L:Cuntz subeq} (ix)
$$
d_\tau\Big(\big( \sum_n^\infty b_k-\delta\big)_+\Big)\le d_\tau\Big(\big( \underset{ k\ge n, \, k\text{ even}}{ \sum} b_k-\frac{\delta}{2}\big)_+\Big)+d_\tau\Big(\big(\underset{k\ge n, \, k\text{ odd}}{ \sum} b_k-\frac{\delta}{2}\big)_+\Big).$$
The series of the even and odd terms separately are diagonal and  dominated by $B$, hence  they still belong to $\Icon$. Thus it is enough to assume that $\sum_{k=1}^\infty b_k$ itself is diagonal. 

Then $(B-\frac{\delta}{2})_+\in K_c$  by  Proposition \ref {P: Kc is hered cone} (ii), hence $$\widehat{(B-\frac{\delta}{2})_+}= \sum_{k=1}^\infty \widehat{(b_k-\frac{\delta}{2})_+}\in \Aff_+.$$ Since also  $\widehat{(b_k-\frac{\delta}{2})_+}\in \Aff_+$ for every $k$, by Dini's Theorem this series converges uniformly.
But then
$$d_\tau\left(\left( \sum_n^\infty b_k-\delta\right)_+\right)= \sum_n^\infty d_\tau((b_k-\delta)_+)\le \frac{2}{\delta}  \sum_n^\infty \tau((b_k-\frac{\delta}{2})_+) \to 0$$
uniformly on $\TA$. 
\item [(ii)] Again, by Lemma  \ref {L:Cuntz subeq} (ix), for every  $0< \delta'< \delta$ we have
$$d_\tau\left(\left( \sum_{k=1}^\infty b_k-\delta\right)_+\right)\le d_\tau\left(\left( \sum_{k=1}^n b_k-\delta'\right)_+\right)+ d_\tau\left(\left( \sum_{n+1}^\infty b_k-(\delta-\delta')\right)_+\right).$$
By (i), we can choose $n$ such that $d_\tau\left(\left( \sum_{n+1}^\infty b_k-(\delta-\delta')\right)_+\right)< \epsilon$ for all $\tau$.
\ep
\bR{R: difference}
If $B\le \|B\|P$ for some projection $P\in \Icon$, as \cite [Lemma 6.4] {KNZCompPos} shows, but  as it is also easy to verify directly,  the uniform convergence in (i) holds also for $\delta = 0$, and hence, (ii) strengthens to the statement that $$d_\tau\left(\left( \sum_{k=1}^n b_k-\delta\right)_+\right)\to d_\tau\left(\left( \sum_{k=1}^\infty b_k-\delta\right)_+\right)\quad\text{uniformly on }\TA.$$ However, these stronger results do not hold in general as it is readily seen by considering $B:=\sum_{k=1}^\infty \frac{1}{k}(e_{k+1}-e_{k})$ for some approximate identity $\{e_n\}$ in a stable algebra $\A$. Indeed then $B\in \A\subset \Icon$, but $d_\tau\left( \sum_{k=n}^\infty b_k\right)=\infty$ for all $n$.
\eR
We are ready now to prove that strict comparison holds for $\Imin$ provided that it holds for $\A$.
\bT{T:Strict comp Imin}
Let $\A$ be a simple, $\sigma$-unital, non-elementary  C*-algebra with strict comparison of positive elements by traces,  $A, B\in (\Imin)_+$ and assume that $B\not \in \A$. If $d_\tau(A)< d_\tau(B)$ for all $\tau\in \TA$ for which $d_\tau(B)< \infty$, then $A\preceq B$.
\eT
\bp
Let $\epsilon >0$. By Theorem \ref {T: Imin=Icont}, $\Imin= \Icon$. Thus by Lemma \ref {L: abs cont in Imin}  there is a $\delta>0$ and $\alpha >0$ such that $$d_\tau\left(\left(A-\epsilon\right)_+\right) +\alpha \le d_\tau\left(\left(B-4\delta\right)_+\right)\quad \forall\, \tau\in \TA.$$ By the assumption that $B\not\in \A$, we can reduce if necessary $\delta$ so that  $(B-4\delta)_+\not \in \A$.
By Theorem \ref {T:bidiag revisited} and Remark \ref {R:bidiag},  $B=\sum_{k=1}^\infty b_k+t$ where $\sum_{k=1}^\infty b_k$ is a strictly converging bi-diagonal series, $t=t^*\in \A$, and $\|t\|< \delta$.
Then by Lemma \ref {L: handy ineq}
\be{e:Cuntz chain} 
(B-4\delta)_+\preceq \left(\sum_{k=1}^\infty b_k-3\delta\right)_+\preceq B\ee
whence by (\ref{e: suub->ineq})  for all $\tau$
$$
d_\tau\left(\left(A-\epsilon\right)_+\right) +\alpha \le d_\tau\left(\left( \sum_{k=1}^\infty b_k-3\delta\right)_+\right).
$$
By Lemma \ref {L: bidiagonal} (ii), there is an $n_1$ such that
\be{e:n1}
d_\tau\left(\left(A-\epsilon\right)_+\right) < d_\tau\left(\left( \sum_{k=1}^{n_1} b_k-2\delta\right)_+\right)~ \forall \tau\in \TA.
\ee
Since $(B-4\delta)_+\not \in \A$,  we have by (\ref {e:Cuntz chain}) that  $\left(\sum_{k=1}^\infty b_k-2\delta\right)_+\not \in \A$.
But then 
\be{e: nonzero}\forall ~n~\exists ~m\ge n \text{ such that }
\left(\sum_n^m b_k -2\delta\right)_+\ne 0.
\ee
Otherwise if there were an $n$ such  that  $\left(\sum_n^m b_k -2\delta\right)_+=0$ for all $m$, the strict convergence $\left(\sum_n^m b_k -2\delta\right)_+\to \left(\sum_n^\infty b_k -2\delta\right)_+$ (see \cite [Lemma 3.1] {KNZCompPos}) would imply that $\left(\sum_n^\infty b_k -2\delta\right)_+=0$, and hence, from Lemma \ref {L:Cuntz subeq} (ix), 
 $$\left(\sum_{k=1}^\infty b_k-2\delta\right)_+\preceq \sum_{k=1}^{n-1} b_k+ \left(\sum_n^\infty b_k -2\delta\right)_+\in \A,$$
 a contradiction. Now starting with the integer $n_1$ just constructed, and by the same argument, define inductively an increasing  sequence of integers $n_k\ge n_{j-1}+2$ such that $$\left(\sum_{n_k+2}^{n_{k+1} }b_k -2\delta\right)_+\ne 0 ~ \forall\, k.$$ Let $d_1:= \sum_{j=1}^{n_1} b_k $ and $d_{k+1}:= \sum_{n_k+2}^{n_{k+1} }b_k $. By construction, $d_nd_m=0$ for $n\ne m$ and 
\be{e:ineq 18}\sum_{k=1}^\infty d_k\le \sum_{j=1}^\infty b_k.\ee
By construction $\left(d_k-2\delta\right)_+\ne 0$ for all $k$ and the function $d_\tau\left(\left(d_k-2\delta\right)_+\right)$ is lower semicontinuous and strictly positive. Let  $$\beta_k :=\min_{\tau\in \TA}d_\tau\left(\left(d_k-2\delta\right)_+\right).$$  By (\ref {e:n1}) we also have 
\be {e:ineq7}d_\tau\left(\left(A-\epsilon\right)_+\right) < d_\tau\big((d_1-2\delta)_+\big) ~ \forall\, \tau.\ee
Now apply Theorem \ref {T:bidiag revisited} and Remark \ref {R:bidiag} to decompose $A$ into the strictly converging sum of a series $\sum_{j=1}^\infty a_k$ and a selfadjoint remainder $a\in \A$ with $a_k\in \A_+$, $a_ka_i=0$ for $|i-j|\ge 2$, and $\|a\|\le \epsilon$. 
By Lemma \ref {L: bidiagonal} (i) we can find a strictly increasing sequence of integers $m_k$ such that
$$
d_\tau \Big(\big(\sum_{j= m_{k}+1}^{\infty}a_k-2\eps\big)_+\Big)<\beta_{k+1} \  \  \ \forall\, \tau\in \TA.
$$
Set $m_o=0$ and $c_k:= \sum_{j= m_{k-1}+1}^{m_{k}}a_k$. 
We claim that
\be{e:strict ineq}
d_\tau\big ((c_k-2\eps)_+\big)< d_\tau\big((d_k-2\delta)_+\big)\quad \forall\, \tau\in \TA, \, k\ge 1.
\ee
For $k=1$ we have
$$
(c_1-2\eps)_+=\left( \sum_{j=1}^{m_1} a_k-2\epsilon\right)_+ \preceq \left( \sum_{j=1}^\infty a_k-2\epsilon\right)_+\preceq \left(A-\epsilon\right)_+.
$$
where the first sub-equivalence follows from Lemma \ref {L:Cuntz subeq} (viii) and the second one from Lemma \ref {L: handy ineq}. Then by (\ref {e: suub->ineq}) and (\ref  {e:ineq7}), $$d_\tau\big ((c_1-2\eps)_+\big)\le d_\tau\left(\left(A-\epsilon\right)_+\right) < d_\tau\big((d_1-2\delta)_+\big),$$
that is, (\ref{e:strict ineq}) holds for $k=1$.
For $k\ge 2$, by Lemma \ref {L:Cuntz subeq} (viii)  and (\ref {e: suub->ineq}) we have for all $ \tau\in \TA$  that
$$
d_\tau\big ((c_k-2\eps)_+\big)\le d_\tau\big ((\sum_{j= m_{k-1}+1}^{\infty}a_k-2\eps)_+\big)<\beta_{k}, 
$$
and hence, (\ref{e:strict ineq}) also holds.

By the strict comparison of positive elements in $\A$, it follows that 
\be{e: comp at A}
(c_k-2\eps)_+\preceq (d_k-2\delta)_+\quad \forall\, k\ge 1.
\ee
By construction, $\sum_{k=1}^\infty a_k=\sum_{k=1}^\infty c_k$ with convergence in the strict topology
and $c_nc_m=0$ for $|n-m|\ge 2$. Thus $C_e:= \sum_{k=1}^\infty c_{2k}$ and $C_o:= \sum_{k=1}^\infty c_{2k-1}$ are two diagonal series also converging strictly and $\sum_{k=1}^\infty a_k= C_e+C_o$. Furthermore, 
$$(C_e-2\eps)_+= \sum_{k=1}^\infty (c_{2k}-2\eps)_+~ \text{and} ~ (C_o-2\eps)_+= \sum_{k=1}^\infty (c_{2k-1}-2\eps)_+.$$
By Proposition \ref {P:diag Cuntz} we have 
\be{e:evenodd}(C_e-3\eps)_+\prec \big( \sum_{k=1}^\infty d_{2k}-\delta)_+~ \text{and} ~ (C_o-3\eps)\prec \big( \sum_{k=1}^\infty d_{2k-1}-\delta)_+.\ee
Therefore
\begin{alignat*}{2}
(A-7\eps)_+&\preceq(C_e+C_o-6\eps)_+ &(\text{by Lemma \ref {L: handy ineq}})\\
&\preceq(C_e-3\eps)_++(C_o-3\eps)_+&(\text{by Lemma \ref {L:Cuntz subeq} (ix)})\\
&\preceq \big( \sum_{k=1}^\infty d_{2k}-\delta)_+\oplus \big( \sum_{k=1}^\infty d_{2k-1}-\delta)_+&(\text{by Lemma \ref {L:Cuntz subeq} (vi)})\\
&= \big( \sum_{k=1}^\infty d_{k}-\delta)_+&(\text{by Lemma \ref {L:Cuntz subeq} (vii)})\\
&\preceq\big( \sum_{k=1}^\infty b_{k}-\delta)_+&(\text{by (\ref {e:ineq 18}), Lemma \ref {L:Cuntz subeq} (viii)})\\
&\preceq B&(\text{by Lemma \ref {L:Cuntz subeq} (iv)})
\end{alignat*}
Since $\eps$ is arbitrary, we conclude that $A\preceq B$.

\ep

\section{An example  where $\Imin \neq \Icon$.}\label{S:Villadsen} From Theorem \ref {T: Imin=Icont}, examples where $\Imin \neq \Icon$ can be found only among ``pathological" algebras that do not have strict comparison of positive elements.  In this section we prove that  the algebras constructed by Villadsen in  \cite{VilladsenSR} provide such examples.  We will largely follow his notations.
Let
$$X_0=\mathbb D^{n_0}\quad\text{and}\quad X_i= X_{i-1}\times  \mathbb C P ^{n_i}~\text{ for } i\in \mathbb N,$$
that is, 
 $$X_i= \mathbb D^{n_0} \times \mathbb C P ^{n_1} \times  \mathbb C P ^{n_2}\times \cdots \times \mathbb C P ^{n_i}.$$
 We will  always assume that 
 \be{e: growth}n_i\ge \sigma(i):=\begin{cases}1&i=0\\i(i!)&i\ge 1\end{cases},\ee 
and hence,  
\be{e: dim}\dim (X_i)=2\sum_{k=0}^i n_k\ge 2 \sum_{k=0}^i \sigma (k)= 2(i+1)!\ee 
This condition, together with the appropriate connecting maps,  will guarantee that the AH  algebra $\A$ constructed in this process will not have {\it slow dimension growth}, which by \cite [Corollary  4.6] {TomsCuntzStab} would imply strict comparison of positive elements.  
 We refer the reader to Villadsen's definition (\cite[pp. 1092-1093] {VilladsenSR}) of the  connecting maps
  \be{e:diag homom}\Phi_{i,i+1}: C(X_i)\otimes\K \to C(X_{i+1})\otimes\K \ee
 and their compositions
 $$ \Phi_{i,j} =  \Phi_{j-1,j}\circ \cdots  \circ\Phi_{i,i+1}: C(X_i)\otimes\K \to C(X_{j})\otimes \K.$$
Identifying as usual projections with complex vector bundles, given  a complex vector bundle $\eta $ over $X_i$, $ \Phi_{i,i+1}(\eta)$ denotes a complex vector bundle over $X_{i+1}$.  Denoting by  $k\eta $ (resp.,  $kq$) the $k$-fold direct sum of the vector bundle $\eta$ (resp., of the projection $q$) with itself, we then have
  \be{e: one step} \Phi_{i,i+1}(\eta)\cong \eta \times \big((i+1)\rank(\eta)\big)\gamma_{n_{i+1}}.\ee
Here $\gamma_{k}$ denotes the universal line bundle over the projective space $\mathbb C P^k$ (see (\ref {e: Euler of universal}) below for a key property of $\gamma_{k}$). Iterating we have for every $j> i$, 
\be{e: many steps}
 \Phi_{i,j}(\eta)\cong \eta \times \frac{\sigma(i+1)}{(i+1)!}\rank(\eta) \gamma_{n_{i+1}}\times  \frac{\sigma(i+2)}{(i+1)!}\rank(\eta) \gamma_{n_{i+2}}  \cdots\times  \frac{\sigma(j)}{(i+1)!} \rank(\eta) \gamma_{n_{j}}.
   \ee
   In particular, since for every $i$ and $j$, $\rank( \gamma_{i})=1$, $ \sum_{k=0}^j\sigma(k)=  (j+1)!$, and    
   $$\rank(\Phi_{i,j}(\eta))= \rank(\eta) \big(1+ \sum_{k=i+1}^j\frac{\sigma(k)}{(i+1)!}\rank(\gamma_{n_k})\big),$$ we then have
   \be{e: rank ij}\rank(\Phi_{i,j}(\eta))= \frac{(j+1)!}{(i+1)!}\rank(\eta) ~ \, \forall j\ge i.\ee 

  Let $\theta$ be a trivial line bundle over $X_0$ and set 
\ba 
&p_i:= \Phi_{0,i}(\theta)\text{ for }i>0;\\
 &\A_i=p_i(C(X_i)\otimes \K)p_i\text{ for }i\ge 0;\\
 &\A= \varinjlim (\A_i, \Phi_{i, i+1}).\end{align*}

Here  $\Phi_{i,i+1}$ denotes the {\it restriction} of  $ \Phi_{i,i+1}$ to $\A_i$. Let  $\Phi_{i,\infty}: \A_i\to \A$ denote the unital embedding  $\A_i  \hookrightarrow \A$. By Villadsen's construction, these maps are injective and we denote by    $\Phi_{\infty,i}: \Phi_{i,\infty}(\A_i) \to \A_i$ the inverse map of $\Phi_{i,\infty}.$ As usual, we will identify $\A_i$ with their images in $\A$ and  focus on the algebraic inductive limit $\bigcup \A_i\subset \A$.

For ease of reference, notice that
\be{e:rank pi} \rank(p_i)= (i+1)!\quad\forall\, i.\ee

By \cite{Dadarlat et all} (see also a short proof in \cite {VilladsenSR}), $\A$ is a simple, unital, AH-algebra and it has a unique tracial state $\tau$. 
Villadsen proved that if  $n_i=n\sigma(i)$ for a fixed $n\in \mathbb N$, then $\A$ has stable rank $n+1$. What interests us here is that by (\ref{e: dim}) and (\ref{e:rank pi}), $\inf_{i}\frac{\dim(X_i)}{\rank(p_i)}\ge 2$ and hence  $\A$ does not have slow dimension growth,  the group $K_0(\A)$ has perforation, and $\A$ does not have strict comparison of projections by its trace. The same holds for other choices of $n_i\ge \sigma (i)$ as readily seen from Villadsen's construction. 

We will show that  $\Imin\neq \Icon$ for the underlying algebra $\St$ and that every element outside $\Icon$ is full if $\sup\frac{n_i}{\sigma(i)}< \infty$ ($\A$ has {\it flat dimension growth}), while this is not the case for an unbounded dimension growth as   $n_i= i!\sigma(i)$. 

To prove these results, we will focus on {\it diagonal} projections of $\M$, i.e. projections of the form $S= \bigoplus_{k=1}^\infty t_ks_k$ where $t_k\in \mathbb N$, $s_k$ is a projection in $\Phi_{k, \infty}(\A_k)$, and $t_ks_k$ is the direct sum of $t_k$ copies of $s_k$. 

To determine if the diagonal projection $S$ is in $ \Icon$ is easy. Since $\A$ has a unique tracial state $\tau$, and hence, $\Icon=I_\tau$, the projection $S$ is in $\Icon$ if and only if $\tau(S)< \infty$, i.e., $\sum_{k=1}^\infty t_k\tau(s_k)< \infty$. 
If $\eta_k= \Phi_ {\infty, k}(s_k) $ is the complex vector bundle over $X_k$ corresponding to $s_k$, i.e.,  $s_k= \Phi_{k, \infty}(\eta_k)$, then by (\ref {e: rank ij}) $ \tau(s_k)= \frac{\rank(\eta_k)}{\rank(p_k)}= \frac{\rank(\eta_k)}{(k+1)!},$
and hence,
\be{e:trace of bundle}\tau(S) = \sum_{k=0}^\infty \frac{t_k\rank(\eta_k)}{(k+1)!}.\ee
To construct a diagonal projection $S\not \in \Imin$ we will make use of algebraic topology tools, more precisely, properties of the Euler classes. For a complex vector  bundle $\eta$ on a compact metric space $X$, $e(\eta)$ will denote the Euler class in the cohomology ring $H^*(X)$. To simplify notations, we will suppress explicit reference to the base space $X$. We start  by recalling that for the universal line bundles $\gamma_{n_{i}}$  used in defining the connecting maps (\ref {e: one step}), we have  \be{e: Euler of universal}
 e(\gamma_{n_{i}})^n\begin {cases} \ne 0& n\le n_i\\ =0& n> n_i.\end{cases}
 \ee

\bL {L: vanishing Euler}
Let $\eta$ be a vector bundle over $X_i$ and let $j>i$.
\item [(i)] If $e(\eta) =0$,  then $e(\Phi_{i, j}(\eta))= 0$.
\item [(ii)] If $e(\eta) \ne 0$ and $\rank(\eta)\le (i+1)!$, then $e(\Phi_{i, j}(\eta))\ne  0$.
\eL
\bp
Recall the fact that the Euler class of  $\Phi_{i,j}(\eta)$ is the cup product of the Euler classes of its components in the Cartesian product in (\ref {e: many steps}) (viewed as vector bundles over $X_j$ via pullbacks of the relevant projection maps).
 That is,  
\ba
e(\Phi_{ij}(\eta))&= e(\eta)e\big(\frac{\sigma(i+1)}{(i+1)!}\rank(\eta)\gamma_{n_{i+1}}\big)\cdots e\big(\frac{\sigma(j)}{(i+1)!}\rank(\eta)\gamma_{n_{j}}\big)\\
&=e(\eta)e(\gamma_{n_{i+1}})^{\frac{\sigma(i+1)}{(i+1)!}\rank(\eta)}\cdots e(\gamma_{n_{j}})^{\frac{\sigma(j)}{(i+1)!}\rank(\eta)}.
\end{align*}
Thus if $e(\eta)$ vanishes, so does $e(\Phi_{ij}(\eta))$. By the Kunneth formula, since the cohomology groups considered have no torsion, it follows that $e(\Phi_{ij}(\eta))\ne 0$ if and only if all the factors in the above decomposition do not vanish. By (\ref {e: Euler of universal}), a necessary and sufficient condition for that to happen is that $n_k \ge \frac{ \sigma(k)}{(i+1)!}\rank(\eta)$  for all $i<k\le j$. By the assumption (\ref{e: growth}), a sufficient condition  is that $\rank(\eta)\le (i+1)!$\ep

Recall that, in each building block $\A_i \otimes \K$, we are identifying projections with vector bundles over $X_i$. Thus if a projection $p$ belongs to $\A_i\otimes \K$ for some $i$  
we associate with it the sequence $\{\eta_j\}_i^\infty$  of the vector bundles  $\eta_j:= \big(\Phi_{ \infty, j}\otimes \id\big)(p)$ over the spaces $X_j$, and $\eta_j= \Phi_{ij}(\eta_i)$ for $j\ge i$.
In view of Lemma \ref {L: vanishing Euler}, it is convenient to set the following definition.
\bD{D: Euler of proj} Let  $p\in  \big(\bigcup _{j=0}^\infty \A_j\big)\otimes \K$ be a projection. We say that
\item [(i)] $e(p)=0$ if $e(\eta_i)= 0$ for some $i$ for which $p\in \A_i\otimes \K$ (and hence $e(\eta_j)= 0$ for every $j\ge i$.)
\item [(ii)] $e(p)\ne 0$ if $e(\eta_j)\ne 0$ for every $j$ for which $p\in \A_j \otimes \K$.
\eD
In order to verify that $e(p)\ne 0$, by Lemma \ref {L: vanishing Euler} it is sufficient  to show that  $e(\eta_i)\ne 0$ and that $\rank( \eta_i)\le (i+1)!$ for the smallest $i$ for which  $p\in \A_i\otimes \K.$ 

\bC{C: subequiv vs Euler}
Let  $q,r\in \big(\bigcup _{j=0}^\infty \A_j\big)\otimes \K$ be projections, $q\preceq r$ and $e(q)=0$. Then $e(r)=0$.
\eC
\bp
There is an $i$ such that $q, r\in \A_i\otimes \K$,   $e\big(\Phi_{\infty, i}\otimes \id)(q)\big)=0$, and the subequivalence $q\preceq r$ holds  within $\A_i\otimes \K$, i.e.,  $r= q'\oplus s$ for some projections  $q', s\in \A_i\otimes \K$ with $q'=vv^*$ and $q=v^*v$ for some $v\in  \A_i\otimes \K$.
But  then
\ba
e\big((\Phi_{\infty, i}\otimes \id)(r)\big)&= e\big((\Phi_{\infty, i}\otimes \id)(q')\big)e\big((\Phi_{\infty, i}\otimes \id)(s)\big)\\&= e\big((\Phi_{\infty, i}\otimes \id)(q)\big)e\big((\Phi_{\infty, i}\otimes \id)(s)\big)=0.
\end{align*}
By Definition \ref {D: Euler of proj}, $e(r)=0$.
\ep

We will construct now two sequence of projections $\{q_i\}$ and $\{r_i\}$ in $\St$ for which $e(q_i)=0$ and $e(r_i)\ne 0$ for all $i$. 

By the definition of $p_i$ it is immediate to find a trivial complex line bundle $\theta_i\le p_i$ over $X_i$. Let $q_i:=\Phi_{i, \infty}(\theta_i)\otimes e_{ii}\in \St$, so that the projections $q_i$ are mutually orthogonal. Then it is clear that 
$$Q:=\bigoplus_{i=1}^\infty q_i\in \M \setminus \St\quad\text{and }\tau(Q)= \sum_{i=1}^\infty \tau( q_i)=  \sum_{i=1}^\infty \frac{1}{(i+1)!}< \infty,$$
and hence, 
\be{e:Q in Icon} Q\in \Icon.\ee Furthermore, by construction, \be{e: Euler qi} e(q_i)=0\quad \forall \, i.\ee 
Next, from the definition of $p_i$ and the construction of the maps $\Phi_{i,i+1}$ in (\ref{e:diag homom}), we see that there is a  complex line bundle  $ \rho_i\in C(X_i)\otimes \K$ with $\rho_i\le p_i $ and $ \rho_i\sim \pi_i^{2*}(\gamma_{n_i})$ where $\pi_i^{2*}$ denotes the pull back map from vector bundles on $\mathbb CP^{n_i}$ to those on the space $X_i$. Thus $e(\rho_i)^k=0$ if and only if $e(\gamma_{n_i})^k=0$, i.e., by (\ref {e: Euler of universal}), if and only if $k> n_i$. When there is no risk of confusion, we will write $\gamma_{n_i}$ for $\rho_i$ as well as for the pullbacks to vector bundles over $X_i$ for $j > i$. Set $$r_i:= \Phi_{i, \infty}(\rho_i)\otimes e_{ii}\in \A_i\otimes \K\subset \St.$$ By definition,  the projections $r_i$ are mutually orthogonal. Set 
$$R:=\sum_{i=1}^\infty r_i$$
 It is the clear that
$R\in \M \setminus \St $, $$\tau(R)= \sum_{i=1}^\infty \tau( r_i)=  \sum_{i=1}^\infty \frac{1}{(i+1)!}< \infty,$$
and hence 
\be {e:R in Icon} R\in \Icon.\ee  

\bL{L:e: Euler nri} $e(nr_i)\ne 0$ for all $n\le \min (n_i, (i+1)!)$. In particular, $e(r_i)\ne 0$ for all $i$.
\eL
\bp
By (\ref {e: Euler of universal}) and the assumption that $n\le n_i$ we have
$e(n\gamma_{n_i})= e(\gamma_{n_i})^n \ne 0$. Moreover, $\rank (n\gamma_{n_i})=n\le (i+1)!$, hence  $e(nr_i)\ne 0$ by Lemma \ref  {L: vanishing Euler} and Definition \ref {D: Euler of proj}.
\ep

\bL{L: rn vs rm} For all integers $j> i$
\item [(i)] $\big(\frac{j!}{i!}\big)r_j\preceq r_i$;
\item [(ii)] $\sum_{k=i+1}^j r_k\preceq r_{i}$;
\item [(iii)] $\sum_{k=i}^j r_k\preceq 2r_{i}$.
\eL
\bp
\item [(i)] By (\ref {e: one step}) we have 
$\Phi_{i,i+1}(\rho_i)\cong \rho_i\times (i+1) \gamma_{n_{i+1}},$
and hence, $(i+1)r_{i+1}\preceq r_i.$ Then (i) follows immediately.
\item [(ii)]  The proof is by induction on $j-i\ge 1$. By (i), $r_{i+1}\preceq (i+1)r_{i+1}\preceq r_i$ so the condition holds for $j-i=1$. Assume  condition (ii) holds for some $j-i>1$ and hence  $\sum_{k=i+2}^{j+1} r_k\preceq r_{i+1}$. Then
$$\sum_{k=i+1}^{j+1} r_k\preceq r_{i+1}\oplus \sum_{k=i+2}^{j+1} r_k\preceq 2r_{i+1}\preceq (i+1)r_{i+1}\preceq r_i$$
 where the last relation in the chain holds by (i).
\item [(iii)] Obvious from (ii).
\ep

\bL{L:rk is order dense}
The sequence $\{r_k\}$ is order dense (see Definition \ref {D:rapidly decreasing}).
\eL
\bp
In view of Lemma \ref {L:Cuntz subeq} (iv) and the density of $\bigcup_{i=1}^\infty\Phi_{i, \infty}(\A_i)$ in $\A$, in order to show that  $\{r_k\}$ is order dense in $\St$, it is enough to show that for every $i, h\in \mathbb N$  and $0\ne a\in p_i(C(X_i)\otimes M_h(\mathbb C))_+p_i$ there is some $j>i$ such that $r_j\preceq \Phi_{i,j}(a)$. 

To do that we need to examine more closely the construction of the connecting maps $\Phi_{i,i+1}$ and their iterations $\Phi_{i,j}$. We again refer the reader to the definition in \cite {VilladsenSR} and also to  \cite {KuchNg}. Disregarding the isomorphism between $\K\otimes \K$ and $\K$, we may view $\Phi_{i,i+1}(a)$ to be in the following matrix form:

$$\begin{pmatrix} a\circ \pi_{i+1, i}&&&\\
&a(\pi_{i+1, i}(y_{i+1}^1))\otimes r_{i+1}^1&&\\
&&&\ddots&\\
&&&&a(\pi_{i+1, i}(y_{i+1}^{i+1}))\otimes r_{i+1}^{i+1}\end{pmatrix}$$
where $r_{i+1}^k$ are mutually orthogonal projections all equivalent to $ r_{i+1}$, $\pi_{i+1, i}$ denotes the projection from $X_{i+1}$ onto $X_i$, and the points $y_{j}^k\in X_j$ are chosen so that the collection of their projections $\{\pi_{j, i}(y_j^k)\mid 1\le k\le j,\, j\ge i\}$ is dense in $X_i$ for every $i$. Since $a$ is a continuous, there is a $j>i$ and a $1\le k\le j$ such that $a(\pi_{j, i}(y_j^k))\ne 0$. But then, 
$0\le a(\pi_{j, i}(y_j^k))\otimes r_j^k\le \Phi_{i,i+1}(a)$. By diagonalizing $a(\pi_{j, i}(y_j^k))$, we can find a $\lambda >0$ and a rank one projection $s$ such that $\lambda s\otimes r_{j}^k\le \Phi_{i,j}(a)$, and hence, $r_j\preceq \Phi_{i,j}(a)$. This proves the claim.

\ep
\bC{C: R generates Imin}
The projection $R$ belongs to $\Imin$, and hence, it generates $\Imin$.
\eC
\bp
Let $e_n:= 1_\A\otimes \sum_{k=1}^n e_{kk}$, then $e_n$ is an approximate identity of $\St$. By Lemma \ref {L:rk is order dense}, the sequence $\{r_k\}$ is order dense, and hence,  by Lemma \ref {L: rn vs rm} it is thin. But then $R\in K_o(\{e_n\})\subset \Imin$  by Lemma \ref {L: diag}. Since $R\not \in \St$ and $\Imin$ is minimal among the ideals properly containing $\St$, it follows that $R$ generates $\Imin$.
\ep

\bT{T: Imin not Icont}
The projection $Q$ does not belong to $\Imin$, and hence, $\Imin \ne \Icon$.
 \eT
 \bp
 Assume by contradiction that $Q\in \Imin$. By Corollary \ref {C: R generates Imin}, $\Imin= I(R)$, and hence there is an $n\in \mathbb N$ such that $Q\le n R$, i.e.,  $\bigoplus_{k=1}^\infty q_k\preceq \bigoplus_{k=1}^\infty nr_k$. 
Choose $i$ such that $n\le \sigma (i-1)$. Then  $n\le \min(n_{i-1}, i!)$ by the assumption (\ref{e: growth}) and hence $e(nr_{i-1})\ne 0$ by Lemma \ref{L:e: Euler nri}.
On the other hand, by Proposition \ref {P: from subeq} there are $m, j\in \mathbb N$, $j\ge i$, such that $q_m\preceq \bigoplus_{k=i}^j nr_k$. By Lemma \ref {L: rn vs rm} (ii),  $q_m\preceq n r_{i-1}$ and since 
$e(q_m)=0$,  it follows from Corollary \ref {C: subequiv vs Euler} that $e(nr_{i-1})=0$, a contradiction.
 \ep 
\bR{R:not strict comp}
\item [(i)] A consequence of Lemma \ref {L:rk is order dense} is the known fact that Villadsen's algebras have the (SP) property (e.g., see the proof of the (SP) property for the Villadsen's type algebras studied in \cite {RordamRR}.) 
\item [(ii)]
The same argument in the proof of Theorem \ref {T: Imin not Icont} shows that $q_m\not \preceq r_i$ for every $m, i\in \mathbb N$ which is an illustration of the well known fact that strict comparison of projections does not hold in $\St$. \eR 

Notice that so far we have only assumed that $n_i\ge \sigma (i)$. We can obtain more if we assume that $\A$ has {\it flat dimension growth}, that is $\sup \frac{\dim(X_i)}{\rank(p_i)}< \infty$, (see \cite [Definition 1.2] {TomsFlat}), which are exactly  Villadsen's {\it finite} stable rank algebras studied in \cite{VilladsenSR}.

\bT{T: Icon maximal} 
Assume that $\A$ has flat dimension  growth, then $\Icon$ is the largest proper ideal of $\M$.
\eT
\bp
To prove that  $\Icon$ contains every proper ideal of $\Ma$, it suffices to prove that  if   $S\in \Ma_+\setminus \Icon$ then $S$ is full, namely $I(S)=\Ma$.   Assume without loss of generality that $\|S\|=1$. By Theorem \ref {T:bidiag revisited} and Remark \ref {R:bidiag}, $S=D_e+D_o+a$ where $a=a^*\in \A\subset \Icon$ and $D_e$ and $D_o$ are diagonal series. Then at least one of the two series must also not belong to $\Icon$.  To simplify notations, assume that $S$ itself is diagonal, namely $S=\bigoplus_{k=1}^\infty s_k$ where $s_k\in \St_+$ for every $k$ and the series converges strictly. Furthermore,  find $\delta>0$ for which $\tau(S-\delta)_+= \infty$. 
Let $M:=\sup \frac{\dim(X_i)}{\rank(p_i)}$ and choose an increasing subsequence $m_k$ such that $\sum_{j=m_k+1}^{m_{k+1}}\tau((s_j-\delta)_+) > \frac{M}{2}+2$. To simplify notations, assume $m_k=k$, i.e., $\tau((s_k-\delta)_+) > \frac{M}{2}+2$ for every $k$.
It was proven in \cite[Lemma 2.5]{KuchNg} that for every $0\ne c\in (\A\otimes M_n)_+$ and $\eps>0$, there is a projection $q$ with $|\tau(q)-\underset{k\to\infty}{\lim }\tau(c^{1/k}) |< \eps$ and $q\in \overline{c(\A\otimes M_n)c}$, and hence, $q\preceq c$. While  the standard assumption in \cite{KuchNg} was that $n_i=\sigma(i)$, no conditions on $n_i$ were used in the proof of that lemma. Moreover, it is routine to extend that lemma to $0\ne c\in (\A\otimes \K)_+$. Thus we can
find projections $q_k\preceq (s_k-\delta)_+$, such that for all $k$ $$\tau(q_k)> d_\tau((s_k-\delta)_+)- \frac{1}{2^k} \ge \tau((s_k-\delta)_+))- \frac{1}{2^k}> \frac{M}{2}+\tau(1_\A\otimes e_{kk}).$$ 
By \cite [Definition 2.1]{TomsFlat}, $M\ge \drr(\A)$  (we refer the reader to \cite{TomsFlat} for the definition of the {\it dimension-rank ratio} of $\A$)  and by \cite [Theorem 3.10] {TomsFlat} it follows that $$1_\A\otimes e_{kk}\preceq q_k\preceq (s_k-\delta)_+\quad \forall\, k.$$
Then  by Proposition \ref {P:diag Cuntz} we have that
$$1_{\Ma}= \bigoplus _{k=1}^\infty 1_\A\otimes e_{kk}\preceq \bigoplus _{k=1}^\infty s_k=S$$
       which proves that $S$ is full.                       
\ep

 Without the flat dimension growth condition, the conclusion of Theorem \ref {T: Icon maximal} no longer necessarily holds. 
To show that, we first we need the following refinement of Lemma \ref {L: rn vs rm}.

 \bL{L:finite sum}
  Let $\eta= \bigoplus _{k=i}^j\Phi_{k, j}(t_k\gamma_{n_k})$, where $i\le j$ are integers and $t_k$ is a monotone nondecreasing sequence of integers. For every $j'\ge j$ we have $$\Phi_{j,j'}(\eta)\cong m_i \gamma_{n_i}\times  m_{i+1} \gamma_{n_{i+1}}\times \cdots m_{j'} \gamma_{n_{j'}}$$ where $m_k\in \mathbb N$ and 
  $$m_k\le \begin{cases}t_i&k=i\\t_k\big(1+ e\frac{\sigma(k)}{(i+1)!}\big)&i+1\le  k\le j\\
 t_je\frac{\sigma(k)}{(i+1)!}&j+1\le  k\le j' .\end{cases}$$
     \eL 
\bp
From (\ref {e: many steps}) we have for every $j'\ge j$
$$\begin{array}{crrrrrrrrrrrrrrrr}
\Phi_{i, j'}(t_i\gamma_{n_i})&\cong& t_i\gamma_{n_i}&\times& t_i\frac{\sigma(i+1)}{(i+1)!}\gamma_{n_{i+1}} &  \times & \cdots & \times &  t_i \frac{\sigma(j)}{(i+1)!}\gamma_{n_{j}} & \times & \cdots & \times & t_{i}\frac{\sigma(j')}{(i+1)!}\gamma_{n_{j'}} \\
\Phi_{i+1, j'}(t_{i+1}\gamma_{n_{i+1}})&\cong& & & t_{i+1}\gamma_{n_{i+1}}  & \times & \cdots & \times & t_{i+1} \frac{\sigma(j)}{(i+2)!}\gamma_{n_{j}} & \times & \cdots & \times &t_{i+1}\frac{\sigma(j')}{(i+2)!} \gamma_{n_{j'}}\\
\Phi_{i+2, j'}(t_{i+2}\gamma_{n_{i+2}})&\cong& & &    &  & \cdots &\times &~ t_{i+2} \frac{\sigma(j)}{(i+3)!}\gamma_{n_{j}}&\times&\cdots&\times  &t_{i+2}\frac{\sigma(j')}{(i+3)!} \gamma_{n_{j'}}\\
\vdots&& &   \vdots &  & & &    \vdots&&&&\vdots \\
\Phi_{j, j'}(t_{j}\gamma_{n_{j}})&\cong& & &  &  &  && t_{j} \gamma_{n_{j}}&\times&\cdots&\times &t_{j}\frac{\sigma(j')}{(j+1)!} \gamma_{n_{j'}}
\end{array}$$
Recall that if $\rho_1\cong s_1\alpha\times t_1\beta$ and $\rho_2\cong s_2\alpha\times t_2\beta$ for some complex vector bundles $\alpha$ and $\beta$ on spaces $X$ and $Y$, and integers $s_1$, $s_2$, $t_1$, $t_2$, then $$\rho_1\oplus \rho_2\cong (s_1+s_2)\alpha \times (t_1+t_2)\beta.$$
Thus by summing the integer multipliers of the universal bundles $\gamma_{n_k}$ we obtain that
$$m_k=\begin{cases}t_i&k=i\\t_k+ \sigma(k)\sum_{h=i+1}^{k}\frac{t_{h-1}}{h!}&i+1\le  k\le j\\
 \sigma(k)\sum_{h=i+1}^{j}\frac{t_{h-1}}{h!}&j+1\le  k\le j'
\end{cases}$$
By using the Lagrange remainder of the Taylor series for the exponential function, we see that  
$\sum_{h=i+1}^{k}\frac{1}{h!}\le \frac{e}{(i+1)!}$. This inequality together with the monotonicity of the sequence $t_k$ establishes the claim. \ep  

\bP{P: Rinfty}
Let $R_\infty:= \bigoplus _{k=1}^\infty k! r_k$. Then  $R_\infty \not \in \Icon$. If  $ n_k\ge k!\sigma(k)$, then $Q\not \in I(R_\infty).$
\eP
\bp
Clearly $R_\infty\in \M\setminus \St$ is a projection and $R_\infty \not \in \Icon$ follows from $\tau(R_\infty)= \sum_{k=1}^\infty \frac{k!}{(k+1)!}=\infty$. To show that $Q\not \in I(R_\infty)$ we reason as in the  proof of Theorem \ref {T: Imin not Icont}. For every $n\in \mathbb N$, choose $i$ such that $(i+1)!\ge 2en$ and let $j\ge i$.  Let $\eta$ be the complex vector bundle over $X_j$ corresponding to $\sum_{k=i}^j nk!r_k$. 

Then $\eta \cong \bigoplus _{k=i}^j\Phi_{k, j}(nk!\gamma_{n_k})$, and hence, by Lemma \ref {L:finite sum}, $$\Phi_{j,j'}(\eta) \cong nm_i \gamma_{n_i}\times  nm_{i+1} \gamma_{n_{i+1}}\times \cdots n m_{j'} \gamma_{n_{j'}}$$
Since $$nm_k\le \begin{cases}ni! &\text{for }~k=i\\\
 n k!(1 + e \frac{\sigma(k)}{(i+1)!})&\text{for }~i<k\le j\\
 \frac{en}{(i+1)!} j! \sigma(k)& \text{for }~j+1<k\le j'\end{cases}\le k!\sigma(k)\le  n_k.$$
we see that $e(\Phi_{j,j'}(\eta))\ne 0$ for every $j'\ge j$. 
Thus  $e(n\bigoplus_{k=i}^j t_kr_k)\ne 0$ by Definition \ref {D: Euler of proj}. Reasoning as in the proof of Theorem \ref {T: Imin not Icont}, we then conclude that $Q\not \in I(R_\infty)$.\ep

\end{document}